\documentclass[11pt]{article}
\usepackage{authblk}
\usepackage{amssymb,latexsym,amsmath,graphicx,amsfonts,amsthm,xy,eucal,mathrsfs,fullpage, epsfig}
\usepackage{amscd,subfigure,color,setspace,cite,mathtools}
\usepackage[font=footnotesize]{caption}
\usepackage{float}
\usepackage{verbatim}
\usepackage[running]{lineno}
\usepackage{tikz}
\usetikzlibrary{mindmap,trees,graphs,shadows,backgrounds,shapes.arrows,calc,positioning}
\usetikzlibrary{arrows.meta}
\usepackage{caption}
\usepackage{epstopdf}
\usepackage[english]{babel}
\usepackage{babel}
\usepackage{fontenc}
\usepackage[colorlinks=true,citecolor=blue]{hyperref}
\usepackage{lineno}
\usepackage{pgf,tikz}
\usepackage{array}
\usepackage{multirow}
\makeatletter
\renewcommand*\env@matrix[1][\arraystretch]{%
  \edef\arraystretch{#1}%
  \hskip -\arraycolsep
  \let\@ifnextchar\new@ifnextchar
  \array{*\c@MaxMatrixCols c}}
\makeatother
\usepackage{float}
\restylefloat{table}

\usetikzlibrary{arrows}
\usetikzlibrary[patterns]
\usetikzlibrary{decorations.pathmorphing}
 \usetikzlibrary{plotmarks}
        \usetikzlibrary{intersections}
        \usetikzlibrary{shadings}
        \usetikzlibrary{calc}
\newlength{\hatchspread}
\newlength{\hatchthickness}
\newlength{\hatchshift}
\newcommand{\hatchcolor}{}
\tikzset{hatchspread/.code={\setlength{\hatchspread}{#1}},
         hatchthickness/.code={\setlength{\hatchthickness}{#1}},
         hatchshift/.code={\setlength{\hatchshift}{#1}},
         hatchcolor/.code={\renewcommand{\hatchcolor}{#1}}}
\tikzset{hatchspread=3pt,
         hatchthickness=0.4pt,
         hatchshift=1pt,
         hatchcolor=black}
\pgfdeclarepatternformonly[\hatchspread,\hatchthickness,\hatchshift,\hatchcolor]
   {custom north west lines}
   {\pgfqpoint{\dimexpr-2\hatchthickness}{\dimexpr-2\hatchthickness}}
   {\pgfqpoint{\dimexpr\hatchspread+2\hatchthickness}{\dimexpr\hatchspread+2\hatchthickness}}
   {\pgfqpoint{\dimexpr\hatchspread}{\dimexpr\hatchspread}}
   {
    \pgfsetlinewidth{\hatchthickness}
    \pgfpathmoveto{\pgfqpoint{0pt}{\dimexpr\hatchspread+\hatchshift}}
    \pgfpathlineto{\pgfqpoint{\dimexpr\hatchspread+0.15pt+\hatchshift}{-0.15pt}}
    \ifdim \hatchshift > 0pt
      \pgfpathmoveto{\pgfqpoint{0pt}{\hatchshift}}
      \pgfpathlineto{\pgfqpoint{\dimexpr0.15pt+\hatchshift}{-0.15pt}}
    \fi
    \pgfsetstrokecolor{\hatchcolor}
    \pgfusepath{stroke}
   }
\usetikzlibrary{arrows}
\usetikzlibrary[patterns]
\numberwithin{equation}{section}
\newtheorem{theorem}{Theorem}[section]
\newtheorem{definition}[theorem]{Definition}
\newtheorem{lemma}[theorem]{Lemma}
\newtheorem{example}[theorem]{Example}
\newtheorem{proposition}[theorem]{Proposition}
\newtheorem{corollary}[theorem]{Corollary}
\newtheorem{remark}{Remark}[section]

\begin{document}
\title{\textbf{Integrating the enveloping technique with the expansion strategy to establish stability }}
\author[1,2]{ Ziyad AlSharawi\thanks{Corresponding author: zsharawi@aus.edu. This work was done while the first author was on sabbatical leave from the American University of Sharjah.}}
\author[2]{ Jose S. C\'anovas}
\affil[1]{\small American University of Sharjah, P. O. Box 26666, University City, Sharjah, UAE}
\affil[2]{\small Universidad Politécnica de Cartagena, Paseode Alfonso XIII 30203, Cartagena, Murcia, Spain}
\date{\today}
\maketitle
\begin{abstract}
In this paper, we focus on finding one-dimensional maps that detect global stability in multidimensional maps. We consider various local and global stability techniques in discrete-time dynamical systems and discuss their advantages and limitations. Specifically, we navigate through the embedding technique, the expansion strategy, the dominance condition technique, and the enveloping technique to establish a unifying approach to global stability. We introduce the concept of strong local asymptotic stability (SLAS), then integrate what we call the expansion strategy with the enveloping technique to develop the enveloping technique for two-dimensional maps, which allows to give novel global stability results. Our results make it possible to verify global stability geometrically for two-dimensional maps.  We provide several illustrative examples to elucidate our concepts, bolster our theory, and demonstrate its application.
\end{abstract}

\section{ Introduction}
Finding a one-dimensional map that can be utilized to prove the global stability in a more complex one-dimensional map is a three-decades-old notion that can be traced back to Cull \cite{Cul2003}, but having a similar technique for multidimensional maps is a morbid notion. Cull developed and formalized the one-dimensional map approach in the so-called ``enveloping technique" \cite{Cul2003,Cul2005, Cul2007, Cul2008}. The notion of enveloping for a one-dimensional map $f$ with a unique fixed point $\bar x$ relies on identifying a map $g$ that shares the same fixed point and satisfies $g(x)>f(x)$ for all $x<\bar x$, while $g(x)<f(x)$ for all $x>\bar x.$  An interesting choice of the map $g$ can be the self-inverse M\"{o}bius transform  $g(x) =\frac{1-ax}{a-(2a-1)x},\; a\in[0,1) $, with the parameter $a$ can be adjusted to fit the specific choice of $f$ \cite{Cul2007,Cul2005}. The strength of the map $g$ lies in its power to force $f$ not to have a $2$-cycle, which is the key to global stability in one dimension. This enveloping technique is used to prove that local stability implies global stability in several one-dimensional biological models. It is important to emphasize that, as shown by a validating example in \cite{Cul2005B}, enveloping by a M\"{o}bius transform is sufficient but not necessary. Within a related notion, it was shown in \cite{Ma-Ma2007} that enveloping does imply the existence of a global Lyapunov function.

Our aim in this paper is to consider $k$-dimensional maps and focus on the existence of a one-dimensional map $g$ that gives an insight into stability analysis. Thus, we consider the $k$-dimensional difference equation
\begin{equation}\label{Eq-F}
x_{n+1}=F(x_n,x_{n-1},\ldots,x_{n-k+1}),\quad k>1,
\end{equation}
where $F\;:\;\mathcal{I}^k\to \mathcal{I}$ is sufficiently smooth and has a fixed point at $\bar x\in \mathcal{I}.$
Here, $\mathcal{I}$ is a metric space that can be $\mathbb{R},$  $\mathbb{R}_+=[0,\infty)$, or a compact interval. The notion of constructing a one-dimensional map $g$ to prove stability in Eq. \eqref{Eq-F} is not new \cite{El-Li2006,El-Lo2008,Lo2010,Al-Al-Am2015}. However, the key issue lies in determining the appropriate approach and discovering the suitable map $g$ that ensures global stability across the widest range of parameters. This is the driving force behind our research here.  Attempts were made in this direction by finding the so-called \emph{dominance condition} \cite{Lo2010,El-Lo2008} or \emph{dominant functions} \cite{Al-Al-Am2015}. Furthermore, the references cited in \cite{El-Li2006,El-Lo2008,Lo2010} provide several explicit and implicit utilizations of one-dimensional maps in stability analysis. El-Morshedy and Lopez established in \cite{El-Lo2008} the fundamental principle for considering one-dimensional maps as a viable option, but the given results lack the mechanism for constructing the map $g.$ However, their technique will be developed here to simplify finding the one-dimensional map $g$.  Given the close relationship between the dominance condition and the enveloping technique in dimension one, we refer to both as the enveloping technique. In \cite{Al-Al-Am2015}, another line of thought was presented. It was based on finding a sequence of one-dimensional maps $\{f_j\}$ that form a non-autonomous system that bounds $F$ and squeezes its orbits within an invariant set that contains the equilibrium.

The smoothness assumption on the map $F$ gives us the freedom to define local asymptotic stability (LAS) based on the spectrum of the Jacobian matrix.  As usual, an equilibrium solution $\{\bar x\}$ of Eq. \eqref{Eq-F} is called hyperbolic if the spectrum of the Jacobian matrix of $F$ at $\bar x$ does not intersect the unit circle. Our interest throughout this paper is limited to the stability of hyperbolic equilibrium solutions, and based on this, an equilibrium solution of Eq. \eqref{Eq-F} is LAS if the spectrum of the Jacobian matrix at the equilibrium solution is contained in the open unit disk $\mathbb{D}$ of the complex plane. This also encompasses the classical one-dimensional case, i.e., when $k=1$ is allowed in Eq. \eqref{Eq-F}. We call an equilibrium solution globally attracting (GA) if all orbits in the domain (or a specified invariant subset of the domain) converge to the equilibrium solution. When LAS and GA are combined, we refer to it as global asymptotic stability (GAS). It is worth mentioning that GA=GAS in the one-dimensional case, but it is not necessarily true in higher dimensions \cite{Se1997}.  $\mathbb{N}$ is used to denote the set of positive integers, while $\mathbb{Z}^+$ denotes the set of non-negative integers, i.e., $\mathbb{Z}^+=\mathbb{N}\cup\{0\}.$

The subsequent sections of this paper are structured as follows: Section Two provides an overview of several stability techniques and their relationship to the expansion strategy and global stability. Section Three links the higher-dimensional system to a one-dimensional map with the same equilibria and local stability. The objective of the one-dimensional map is to utilize it through the enveloping technique to prove global stability, as discussed in Section Four.  In Section Five, we give some examples and applications that show the strength of our approach. Finally, we end the paper by providing a concise conclusion.

\section{Local stability and expansions}
The popular classical approach to analyzing local asymptototic stability in higher-dimensional maps is based on linearization and the eigenvalues of the Jacobian matrix. However, computing the eigenvalues of the Jacobian matrix remains a significant barrier. The popular approach to tackling the eigenvalue problem is the Jury algorithm, which can be cumbersome for large matrices. The authors in \cite{Al-Ca-Ka2024,Al-Ca2024} developed another straightforward computational algorithm based on the embedding technique and the expansion strategy. We extract the idea of the new algorithm from \cite{Al-Ca2024}, and explain it within the context of our paper. Before we embark on that, it is convenient to normalize a nonzero equilibrium point in Eq. \eqref{Eq-F}. So, we 
let $x_n=\bar xy_n,$ and re-write Eq. \eqref{Eq-F} to become
\begin{equation}\label{Eq-F-Normalized}
y_{n+1}=F_0(y_n,\ldots,y_{n-k+1})=\frac{1}{\bar x}F(\bar xy_n,\bar xy_{n-1},\ldots,\bar x y_{n-k+1}).
\end{equation}
Throughout this paper,  we consider $\mathcal{I}=\mathbb{R}_+=[0,\infty)$ and $F_0:\; \mathcal{I}^k\to \mathcal{I} $ (unless specified otherwise).

\subsection{The expansion strategy}
In Eq. \eqref{Eq-F-Normalized}, it is possible to use the same equation and substitute it in place of $y_n.$ This gives us a new equation of larger delay, namely 
\begin{equation}\label{Eq-F1}
\begin{split}
y_{n+1}=&F_0(F_0(y_{n-1},\ldots,y_{n-k}),y_{n-1},\ldots,y_{n-k+1})\\
=&F_1(y_{n-1},\ldots,y_{n-k}).
\end{split}
\end{equation}
When this process is repeated, we obtain (by induction) a sequence of equations 
\begin{equation}\label{Eq-Fj}
y_{n+1}=F_j(y_{n-j},\ldots,y_{n-k-j+1}),\quad j=0,1,\ldots.
\end{equation}
We consider the smoothness of $F_0$ to be inherited in $F_j,$ and consequently, $F_j\;:\;\mathcal{I}^k\to \mathcal{I}.$  Note that if we re-write the scalar equation in Eq. \eqref{Eq-Fj} in vector form, then it will be from $\mathcal{I}^{k+j}$ to $\mathcal{I}^{k+j}$. The notion of expansions was discussed in \cite{Al-Ca2024} and dubbed as the ``\textit{expansion strategy}". This strategy is found to be effective in developing a local stability algorithm that we present here. It is important to note that a fixed point of $F_0$ is also a fixed point of the map $F_j$ for all $j$, but the other way around is not necessarily true \cite{Al-Ca2024}. Additionally, a solution of Eq. \eqref{Eq-F} can also be a solution of Eq. \eqref{Eq-Fj} provided the initial conditions of Eq. \eqref{Eq-Fj} are determined in a specific manner. Let the Jacobian matrix of $F_j$ at a fixed point of $F_0$ be $J_j$ for $j\in \mathbb{Z}^+.$ The relationship between the spectrum of $J_0$ and $J_j$ has been investigated in \cite{Al-Ca2024}, and we utilize that in our local stability analysis.   At the fixed point $\bar y=1,$ define $a_i$ to be the partial derivative of $F_0$ with respect to its $i$th argument, and let $J_0$ be the Jacobian matrix of System \eqref{Eq-F-Normalized}, i.e.,
$$J_0=\left[
      \begin{array}{ccccc}
        a_1 & a_2 & \cdots & a_{k-1} & a_k \\
        1 & 0 & 0 & 0 & 0 \\
        0 & 1 & 0 & 0 & 0 \\
        \vdots & \vdots & \ldots & \vdots & 0 \\
        0 & 0 & 0 & 1 & 0 \\
      \end{array}
    \right].
 $$
Also, let $V_0$ be the column vector identified by the gradient of $F_0$ at the fixed point ($\nabla F_0$), i.e, $V_0=\nabla F_0= \left[\begin{array}{cccc}
a_1 & a_2 & \cdots & a_k \\
\end{array}\right]^t,$
then define $V_{m}=(J_0^t)^mV_0,$ where $t$ represents the transpose. In this case, 
\begin{equation}\label{In-LASCondition}
\|V_m\|_1=\|\nabla F_m(1,\ldots,1)\|_1=\sum_{i=1}^k\left|\frac{\partial}{\partial x_i}F_m(1,\ldots,1)\right|.
\end{equation}
In the sequel, $\|V_m\|_1$ or $\|\nabla F_m\|_1$ will be used to shorten the expressions in Eq. \eqref{In-LASCondition}. 
\begin{theorem}\cite{Al-Ca2024}\label{Th-Al-Ca2024}
Consider Eq. \eqref{Eq-F-Normalized}, and let $F_j$ be the sequence of expansions introduced in Eq. \eqref{Eq-Fj}. If $\{\bar y=1\}$ is a hyperbolic equilibrium solution, then it is LAS for $F_0$ iff $\|\nabla F_m\|_1<1$ for some non-negative integer $m.$
\end{theorem}

The relationship between the condition $\|\nabla F_m\|_1<1$ in Theorem \ref{Th-Al-Ca2024} and the LAS for $F_0$ has been established and clarified in \cite{Al-Ca2024,Al-Ca-Ka2024}. It will also be clear from the illustrative examples we provide in Section \ref{Sec-Enveloping}. However, the function  $F_0(x,y)=-\frac{3}{5}x-\frac{3}{5}y$ (which we consider again in Eq. \eqref{Eq-Linear1}) gives a simple illustrative example. In this case, the zero fixed point is LAS, but $\|\nabla F_0\|_1>1$. On the other hand, the first expansion is $F_1(x,y)=-\frac{6}{25}x+\frac{9}{25}y$, which gives us $\|\nabla F_1\|_1=\frac{3}{5}<1$. We close this subsection by stressing a crucial fact that will be used in the sequel. 
When the gradient vector of $F_m$ at the fixed point satisfies $\|\nabla F_m\|_1<1,$  we say $F_m$ satisfies the LAS condition at the fixed point, and for ease of reference, we formalize the concept in the following definition:
\begin{definition}\label{Def-LASCondition}
Let $U$ be an open subset of $\mathbb{R}^k$ and $F\;:\; U\subseteq\mathbb{R}^k\;\to\mathbb{R}.$ Assume $F$ has a fixed point $\bar x$ in $U.$ if $F\in\mathcal{C}^1(U)$ and $\|\nabla F\|_1<1$ at the fixed point $\bar x$, we say that $F$ satisfies the LAS condition at $\bar x$. Alternatively, we say $\bar x$ is strongly LAS (SLAS). 
\end{definition}

\subsection{Local stability and one-dimensional maps}
Since we are interested in a one-dimensional map that reveals the global stability of Eq. \eqref{Eq-F},  the appropriate one-dimensional map is expected to pass a local stability test before it can be used to assess global stability. Thus, we begin by employing one-dimensional maps and relating local stabilities. A ``naive" attempt can be to consider $g(x)=F(x,x,\ldots,x).$ To clarify this, consider, for instance, $F(x,y)=ax^2e^{-y}$ for some $a>e.$ There are two positive fixed points; one of the fixed points is larger than one, while the other is smaller than one. However, none of them can be hyperbolic and LAS. On the other hand, $g(x)=F(x,x)$ has the same fixed points, but the small positive fixed point is unstable, and the large one is LAS as long as it is less than three. So, selecting a one-dimensional map $g$ that reflects the local stability in $F$ must be well-articulated.

For a square matrix $J=[u_{ij}],$ recall that the maximum row sum and the maximum column sum norms are defined by
$$\|J\|_{\infty} = \max_{1 \leq i \leq k} \sum_{j=1}^{k} |u_{ij}|\quad \text{and}\quad \|J\|_{1} = \max_{1 \leq j \leq k} \sum_{i=1}^{k} |u_{ij}|,$$
respectively.

\begin{proposition}\label{Pr-Rem1}
Consider $F_0$ in Eq. \eqref{Eq-F-Normalized} with a Jacobian matrix $J_0$ at the fixed point $\bar y=1.$ If $\|\nabla F_m\|_1<1$ at $\bar y$ for some $m,$ then $\|J_0^n\|_\infty$ and $\|J_0^n\|_1$ converge to zero as $n\to \infty.$
\end{proposition}
\begin{proof}
From Theorem \ref{Th-Al-Ca2024}, $\|\nabla F_m\|_1<1$ implies $\bar y$ is LAS, and consequently, the spectral radius  of $J_0$ must be less than $1, $ (i.e., $\rho(J_0)<1$). For any induced matrix norm $\|\cdot\|$, Gelfand’s formula implies  
$$\rho(J_0)=\lim_{n\to \infty}\|J_0^n\|^\frac{1}{n}.$$
Now, since the maximum row sum and the maximum column sum norms are induced matrix norms, then $\rho(J_0)<1$ and Gelfand’s formula give us the desired conclusion. 
\end{proof}
Back to Eq. \eqref{In-LASCondition}, having $\|V_m\|_1<1$ for some $m$ and $V_{m}=(J_0^t)^mV_0$ means we can write $(J_0^t)^m=[v_{i,j}]$, and invoke Proposition \ref{Pr-Rem1} to obtain 
\begin{equation}\label{In-Bj}
B_j:=\left|\sum_{i=1}^kv_{i,j}\right|\leq \sum_{i=1}^k|v_{i,j}|\leq \|J_0^m\|_\infty.
\end{equation}
Furthermore,  we can consider $m$ sufficiently large to obtain $B_j$ as small as needed, i.e., $B_j\to 0$.

When $\bar x=1$ is SLAS under $F_m,$ we define the one-dimensional map
\begin{equation}\label{Eq-OneDimensional-g1}
g(t)=F_m\left(\phi_1(t),\phi_2(t),\ldots, \phi_{k}(t)\right),\quad \text{where}\quad \phi_j(t)=\frac{(c_j-b_j+1)+b_jt}{1+c_jt},
\end{equation}
and $c_j,b_j>0$. Observe that $t=1$ is a fixed point of $\phi_j$ for all $j=1,\ldots,k$, and consequently, it is a fixed point of $g.$ Next, denote the partial derivatives of the expanded map $F_m$ by $F_{n,j}.$  Since,
\begin{equation} \label{Eq-DerivativeOf-g}
\begin{split}
g^\prime(1)=&\sum_{j=1}^kF_{m,j}(\phi_1(1),\phi_2(1),\ldots,\phi_k(1))\phi_j^\prime(1)\\
=& \sum_{j=1}^kF_{m,j}(1,1,\ldots,1)\phi_j^\prime(1),
\end{split}
\end{equation}
and $\phi_j^\prime(1)=\frac{b_j-c_j}{1+c_j}$. To achieve $|g^\prime(1)|<1,$ it is necessary to adjust the parameters $b_j$ and $c_j$ to equilibrate the partial derivatives $F_{m,j}$. In the absence of the LAS assumption in $F_0,$ the partial derivatives $F_{m,j}$ may exhibit significant variability, necessitating the use of negative parameters in $\phi_j$ to achieve balance. We avoid this issue and presume LAS in $F_0$ or SLAS in $F_m$.   
 \begin{theorem}\label{Th-Local1}
 Assume the equilibrium $\bar y=1$ of Eq. \eqref{Eq-F-Normalized} is LAS, and consider the maps $\phi_j$ and $g$ as defined in Eqs.  \eqref{Eq-OneDimensional-g1}. The constants $b_j$ and $c_j$  can be chosen so that the equilibrium $\bar y=1$ is LAS for the one-dimensional map $g$.
 \end{theorem}
 \begin{proof}
Suppose that $\bar y=1$ is a LAS equilibrium of Eq. \eqref{Eq-F-Normalized}. By Proposition \ref{Pr-Rem1}, and the consequent discussion,  we consider $B_j$ as defined in Eq. \eqref{In-Bj}.
As in Eq. \eqref{Eq-DerivativeOf-g}, denote the partial derivative of $F_m$ with respect to its $j$th argument by $F_{m,j},$ and investigate $g^\prime(t)$ at $t=1$. We obtain
\begin{align*}
g^\prime(1)=& \sum_{j=1}^ka_j\phi_j^\prime(1)=\sum_{j=1}^ka_j\frac{b_j-c_j}{1+c_j},
\end{align*}
where $a_j=F_{m,j}(1,\ldots,1).$
If $a_j\geq 0,$ we define $b_j=c_j-(1+c_j)B_j,$ and if $a_j<0,$ then we define $b_j=c_j+(1+c_j)B_j.$ In this case, we obtain 
\begin{equation}\label{Eq-ReviewerPoint}
|g^\prime(1)|= \sum_{j=1}^k|a_j|B_j.
\end{equation}
Now, there exists a positive integer $m_0$  that $\bar y=1$ is SLAS under $F_{m_0}$, i.e., $\|\nabla F_{m_0}\|_1<1.$ Also, from Eq. \eqref{In-LASCondition}, we obtain 
$$V_{qm_0}=\left(J_0^t\right)^{(q-1)m_0}V_{m_0}$$
for any positive integer $q$, which gives us 
$$\|V_{qm_0}\|_1\leq \|J_0^{(q-1)m_0}\|_\infty \|V_{m_0}\|_1= \|J_0^{(q-1)m_0}\|_\infty \|\nabla F_{m_0}\|_1.$$
Since we can consider $q$ to be independent of $m_0,$ we consider $m$ sufficiently large in Eq. \eqref{Eq-ReviewerPoint} to obtain 
$$|g(1)|\leq \|J_0^{(q-1)m_0}\|_\infty \|\nabla F_{m_0}\|_1.$$
By Proposition \ref{Pr-Rem1}, we obtain $|g^\prime(1)|\leq \|\nabla F_{m_0}\|_1<1$ as needed.
 \end{proof}
Note that in the proof of Theorem \ref{Th-Local1}, we still have some freedom to choose $c_j,\; j=1,\ldots,k$.  The free parameters can be calibrated to accommodate specific invariant domains. Also, based on Proposition \ref{Pr-Rem1}, we give the following remark:
\begin{remark}\label{Rem-Epsilon}
 By contemplating the choice of $b_j$ in the proof of Theorem \ref{Th-Local1}, and since $B_j$ can be considered arbitrarily small, we can define $b_j=c_j\pm \epsilon,$ i.e., if $a_j\geq 0,$ we define $b_j=c_j-\epsilon,$ and if $a_j<0,$ we define $b_j=c_j+\epsilon.$
\end{remark}

\begin{proposition}
Assume $\bar y$ is LAS in Eq. \eqref{Eq-F-Normalized}, and $b_j$'s are chosen as in the proof of Theorem \ref{Th-Local1}. The one-dimensional map $g$ in Eq. \eqref{Eq-OneDimensional-g1} has $t=1$ as the unique LAS fixed point  iff $g$ satisfies $(g(t)-t)(t-1)<0$ for all $t\neq 1.$
\end{proposition}
\begin{proof}
Theorem \ref{Th-Local1} ensures that $\bar y=1$ is LAS for $g$. Since $g(1)=1$ and $g^\prime(1)<0,$ there exists a neighborhood $\mathcal{N}$ of $1$ such that $g$ is decreasing for all $t\in\mathcal{N}.$ This implies  $(g(t)-t)(t-1)<0$ for all $t\in\mathcal{N}.$  Because $g$ has a unique fixed point, then $\mathcal{N}$ extends to cover the domain of $g.$ Therefore, we obtain $(g(t)-t)(t-1)<0$ for all $t\neq 1.$ To show the converse, observe that $(g(t)-t)(t-1)<0$ implies $g(t)>t$ when $t<1$ and $g(t)<t$ when $t>1.$  Therefore, $t=1$ is the unique fixed point.
\end{proof}
\begin{proposition}
If the map $F$ in Eq. \eqref{Eq-F} is monotonic in each of its arguments, and $b_j$'s are chosen as in the proof of Theorem \ref{Th-Local1}, then the one-dimensional map $g$ in Eq. \eqref{Eq-OneDimensional-g1} is decreasing with a unique fixed point.
\end{proposition}
\begin{proof}
 Because
$$g^\prime(t)=\sum_{j=1}^kF_{m,j}(\phi_1(t),\phi_2(t),\ldots,\phi_k(t))\phi_j^\prime(t),$$
 then $g^\prime$ is negative at all fixed points. However, since $g$ is continuous, it must increase again before it can have another fixed point. Hence, $g$ is decreasing and the fixed point must be unique.
\end{proof}
 We illustrate the ideas of this section in the following example.
\begin{example}\label{Ex-1}\rm
Consider $x_{n+1}=F_0(x_n,x_{n-1})=x_ne^{b(1-x_{n-1})},$ where $b>0.$ The fixed points of $F_0$ are zero and $\bar x=1.$ Based on the above approach, we have $a_1=1,a_2=-b$ at $\bar x$ and, consequently,
we obtain
$$J=\left[
      \begin{array}{cc}
       1 &-b\\
        1&0\\
      \end{array}\right],
      \quad V_0=\left[
      \begin{array}{c}
        1 \\
        -b\\
      \end{array}\right],
\quad V_1=\left[
      \begin{array}{c}
        1-b \\
        -b\\
      \end{array}\right],
\quad V_2=\left[
      \begin{array}{c}
        1-2b \\
        -b(1-b)\\
      \end{array}
    \right].$$
Therefore, 
$\|\nabla F_0\|_1=\|V_0\|_1>1,$ $\|\nabla F_1\|_1=\|V_1\|_1=1$ and $\|\nabla F_2\|_1=\|V_2\|_1<1$ when $0<b<1.$ So, the positive fixed point is LAS when $0<b<1.$ Now, we have $B_1=|1-2b|$ and $B_2=|1-b|.$ This gives us $b_1=c_1-(1+c_1)|1-2b|$ and $b_2=c_2+(1+c_2)|1-b|.$ Therefore, we obtain
\begin{align*}
\phi_1(t)=& \frac{b_1(t-1)+c_1+1}{1+c_1t }\\
=&\frac{\left(c_1-(1+c_1)|1-2b|\right)(t-1)+c_1+1}{1+c_1t },\\
\phi_2(t)=&\frac{b_2(t-1)+c_2+1}{1+c_2t }\\
=& \frac{\left(c_2+(1+c_2)|1-b|\right)(t-1)+c_2+1}{1+c_2t },
\end{align*}
and consequently,
$$g(t)=\phi_1(t)\exp(b(1-\phi_2(t)))\quad\text{where}\quad g(1)=1.$$
Furthermore, it is straightforward to check that
$$g^\prime(1)=-b|b-1|-|2b-1|<0.$$
\end{example}

 Having the LAS condition satisfied as in Definition \ref{Def-LASCondition} has interesting and favorable outcomes. For instance, a fixed point of the multi-dimensional map $F_0$ could be LAS, while it is unstable for the one-dimensional map $g(x)=F_0( x,\ldots, x)$. This can be made clear by the simple example 
\begin{equation}\label{Eq-Linear1}
    x_{n+1}=F(x_n,x_{n-1})=\frac{3}{5}x_n-\frac{3}{5}x_{n-1}.
\end{equation}
However, forcing the LAS condition on $F$ at a fixed point forces $g(x)=F( x,\ldots, x)$ to be LAS at the fixed point. For ease of reference, we confirm this fact in the following proposition.
\begin{proposition}
Consider the map $F_0$ in Eq. \eqref{Eq-F-Normalized}. If $\bar y=1$ is SLAS under $F_0,$ then it is LAS under the one-dimensional map $g(x)=F_0(x,\ldots,x).$
\end{proposition}
\begin{proof}
This follows right from the triangle inequality, i.e.,
$$|g^\prime(1)|=\left|\sum_{i=1}^k\frac{\partial}{\partial x_i}F_0(1,\ldots,1)\right|\leq \|\nabla F_0(1,\ldots,1)\|_1.$$
\end{proof}
In conclusion,  the dynamics of the one-dimensional map $g(x)=F( x,\ldots, x)$  can be more relevant (at least locally) to the dynamics of $F$ under the LAS condition. 
\section{ Global stability techniques}

Unlike LAS, there is currently no all-encompassing solution for tackling global stability. Varying situations require distinct techniques such as Lyapunov functions, embedding, or enveloping. In this section, we focus on the enveloping technique and integrate it with the expansion strategy, then align them with the goal of our paper. We begin with a mind map that illustrates the interconnection between the ideas:  

\subsection{ A mind map }
To present a cohesive analysis incorporating several methodologies, it is beneficial to envision a mind map that delineates the connections.
\definecolor{DCO}{RGB}{255,131,211}
\begin{figure}[htbp]
\centering
\begin{tikzpicture}[thick, scale=0.7,
 every node/.style={scale=0.8},
 outer sep=2pt
  ]
  \path[mindmap,concept color=purple!10,text=black]
    node[concept] (LS) {\bf Local Stability}
    [clockwise from=0]
    child[concept color=white,level distance=6cm] {
      node[concept] (GS) {\textcolor{red}{\bf Global Stability}}
      [clockwise from=0]
      child [concept color=orange!20!green!10,text=black,level distance=6cm,font=\fontsize{8pt}{10pt}\selectfont] {node[concept,scale=1.4] (ExT) {{\bf  Expansion Strategy}}}
      }
    child [concept color=orange!20, grow=30, level distance=7.5cm, text=black]
          {node[concept] (EnT) {\bf Enveloping Technique}}
    child [concept color=orange!20, grow=330, level distance=7.5cm, text=black]
          {node[concept] (EmT) {\bf Embedding Technique}};
\draw [blue, arrows = {-Stealth}] (EnT) -- (GS);
\draw [blue, arrows = {-Stealth}] (EnT) -- (LS);
\draw [blue, arrows = {-Stealth}] (EmT) -- (GS);
\draw [blue, arrows = {-Stealth}] (EmT) -- (LS);
\draw [blue, arrows = {-Stealth}] (GS) -- (LS);
\draw [blue, arrows = {-Stealth}] (EnT) -- (ExT);
\draw [blue, arrows = {-Stealth}] (ExT) -- (EnT);
\draw [blue, arrows = {-Stealth}] (EmT) -- (ExT);
\draw [blue, arrows = {-Stealth}] (ExT) -- (EmT);
\begin{scope}[on background layer, every node/.style={scale=.85}]
  \path[mindmap,concept color=black,text=black]
    node[concept] {\bf Root}
    [clockwise from=0]
    child[concept color=blue,level distance=6cm] {
      node[concept] {\textcolor{white}{\bf Global Stability}}
      [clockwise from=0]
      child [concept color=red,text=black,level distance=6cm,font=\fontsize{8pt}{10pt}\selectfont] {node[concept,scale=1.4] {{\bf Expansion Strategy}}}
     }
      child [concept color=red, grow=30, level distance=7.5cm, text=black]
            {node[concept] {\bf Enveloping Technique}}
      child [concept color=red, grow=330, level distance=7.5cm, text=black]
           {node[concept] {\bf Embedding Technique}};
\end{scope}
\begin{pgfonlayer}{background}
      \path (EmT) to [circle connection bar switch color=from (orange!30) to (orange!20)] (GS);
       \path (EmT) to [circle connection bar switch color=from (orange!30) to (orange!20!green!30)] (ExT);
       \path (EnT) to [circle connection bar switch color=from (orange!30) to (orange!20!green!30)] (ExT);
        \path (EnT) to [circle connection bar switch color=from (orange!30) to (yellow!30)] (GS);
\end{pgfonlayer}
\end{tikzpicture}
\caption{This diagram is a cognitive map showcasing the ideas presented in this paper. In a prior study \cite{Al-Ca2024}, the authors combined the expansion strategy and the embedding technique to achieve stability outcomes. This paper aims to complete the diagram, integrate the expansion strategy with the enveloping technique, and achieve new stability outcomes.}\label{Fig-MindMap}
\end{figure}
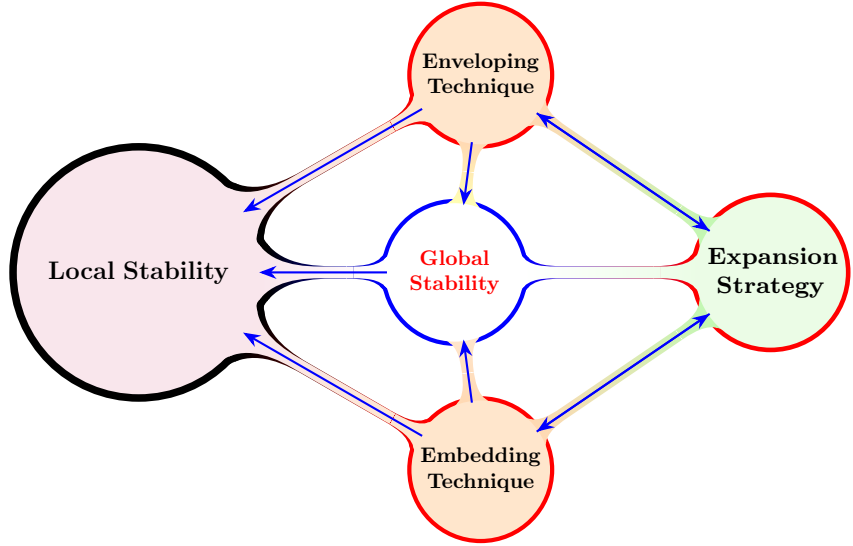

At this juncture, we provide a concise overview of the correlation between global stability, embedding technique, and the expansion strategy. When the map $F$ in Eq. \eqref{Eq-F} is monotonic in each of its arguments, the dynamical system can be embedded in a higher-dimensional monotonic system with respect to a specific partial order. This method works well for achieving global stability as long as $F$ has an invariant compact box and the embedded system does not create any pseudo-fixed points of $F$ \cite{Al-Ca-Ka2024B,Sm2008}. When one of the conditions is invalid, it is possible to increase the delay via successive substitutions and obtain a higher-dimensional system \cite{Al-Ca2024}. This strategy has the potential to generate a higher-dimensional system that enhances the advantages of the embedding technique. The integration of the expansion strategy with the embedding technique was also used in \cite{Al-Ca-Ka2024,Al-Ca2024} to address local stability when the spectrum of the Jacobian matrix is not easy to compute.

Next, we embark on the main task of building the top part of the graph in Fig. \ref{Fig-MindMap}. As mentioned in the introduction, the enveloping technique was initiated for global stability in one-dimensional maps and generalized to higher-dimensional maps by El-Morshedy and Jim\'enez-L\'opez in \cite{El-Lo2008,Lo2010}. The idea about their proposed one-dimensional map $g$ is based on the following:  Let $x_1,\ldots x_k$ be the initial conditions in Eq. \eqref{Eq-F} and,  define $\alpha:=\min\{x_1,\ldots,x_k\}$ and $\beta:=\max\{x_1,\ldots,x_k\}.$ The suitable candidate $g$ must have the same fixed point $\bar x$ and satisfy $(g(x)-x)(x-\bar x)<0$ for all $x\neq \bar x.$ Also, there exists $x\in[\alpha,\beta]$ such that if $F(x_1,\ldots,x_k)\geq \beta,$ then $g(x)\geq F(x_1,\ldots,x_{k}),$ while if $F(x_1,\ldots,x_k)\leq \alpha,$ then $g(x)\leq F(x_1,\ldots,x_{k}).$ Roughly speaking, whenever an element of the orbit of $F$ escapes  the interval $[\alpha,\beta],$ the one-dimensional map $g$ guarantees that it does not escape too far. Several illustrative cases were presented in \cite{El-Lo2008} to show the feasibility of identifying a suitable map $g$.

\subsection{One-dimensional maps through enveloping}
It is well known that a unique fixed point of a continuous map $f$ on a bounded interval $I$ with no $2$-cycles must be GAS \cite{Co1955,Po-Ro1982}. However, we are interested in a slightly different or tailored version here.  We need to drop the condition of having a bounded interval and, instead, we assume the unique fixed point is LAS.  The following proposition addresses this issue. 
\begin{proposition}\label{Pr-2CycleGAS}
Let $f\;:\; [0,\infty)\to[0,\infty)$ be continuous, and assume it has a unique fixed point $\bar x.$ If $f$ has no $2$-cycle, then local stability implies global stability. 
\end{proposition}
\begin{proof}
Because $f$ has a unique fixed point and it is LAS, then the slope of $f$ at $\bar x$ is between $-1$ and $1$, and $f$ crosses $y=x$ at the unique fixed point $\bar x$. Thus, we must have
\begin{equation}\label{In-IneqStability}
(f(x)-x)(x-\bar x)<0\quad\text{for all}\quad x\neq \bar x.
\end{equation}
Define $M:=\max\{f(x):\; 0\leq x\leq \bar x\},$ and let $I=[0,M].$ We have, $M\geq \bar x.$ First, we show that $f:\;I\;\to I.$ Let $x_0\in I.$ If $x_0\leq \bar x,$ then clearly $f(x_0)\leq M$ and $f(x_0)\in I.$ So, consider $x_0>\bar x.$ In this case, Inequality \ref{In-IneqStability} makes $f(x_0)<x_0\leq M$, and thus $f(x_0)\in I.$ Because $f$ maps a bounded interval into itself, and $f^2=f\circ f$ has a unique fixed point, then $I$ is contained in the basin of attraction of $\bar x.$

Next, let $x_0>M.$ If $f^j(x_0)\leq M$ for some $j,$ then $x_0$ also  belongs to the basin of attraction of $\bar x.$ On the other hand, if $f^j(x_0)>M$ for all $j$,  then from Inequality \eqref{In-IneqStability}, we obtain 
$$f^{n+1}(x_0)<f^n(x_0)<\ldots<x_0.$$
This gives a decreasing sequence that must converge to the unique fixed point, $\bar x$. Hence, all orbits must converge to $\bar x,$ and consequently, $\bar x$ is GAS.  
\end{proof}
Note that the condition of local stability in Proposition \ref{Pr-2CycleGAS} is crucial because $f(x)=\alpha x,\; \alpha>1$ maps $[0,\infty)$ onto itself, but the unique fixed point is a repeller. Since proving the existence or non-existence of a $2$-cycle can be challenging, one can attempt to envelop the one-dimensional map with another more straightforward one-dimensional map that achieves the GAS task.  
\begin{theorem}\cite{Cul2008}\label{Th-Cul2008}
Suppose $f\;:\; [0,\infty)\to [0,\infty)$ is continuous and have a unique positive fixed point at $\bar x=1.$ If $g$ is a self-inverse map that satisfies $g(x)>f(x)$ for $0<x<1$ and $g(x)<f(x)$ for $x>1,$ then  $\bar x=1$ is a GAS for $f.$  
\end{theorem}
The authors in \cite{El-Li2006,El-Lo2008,Lo2010} developed the idea of dominance for multi-dimensional maps, which is a generalization of the notion of enveloping. However, the given characterization of the dominant function lacks a mechanism to find a suitable map. We first cite a version of the main result in \cite{El-Lo2008,Lo2010}, then give our more stringent definition of the enveloping concept to establish well-defined mechanisms for finding a suitable enveloping map. 

 \begin{theorem}\cite{El-Lo2008}\label{Th-Lopez}
Let $\mathcal{I}$ be a nonempty interval of $\mathbb{R}$ and $\overline{\mathcal{I}}$ be its closure. Consider Eq. \eqref{Eq-F-Normalized} in which $F_0:\mathcal{I}^k\to \mathcal{I},$ and  let 
$\alpha:=\min\{y_0,\ldots,y_{-k+1}\},$ $ \beta:=\max\{y_0,\ldots,y_{-k+1}\}.$
Suppose there exists a continuous map $g\;:\;\mathcal{I}\to \overline{\mathcal{I}}$ that shares the fixed point $\bar y=1.$  If there exists $t\in[\alpha,\beta]$ such that $g(t)\leq y_1<\alpha$ or $g(t)\geq y_1>\beta$ whenever $y_1=F_0(y_0,y_{-1},\ldots,y_{-k+1})\not\in[\alpha,\beta],$ then the global attractivity of $\bar y$ for $g$ implies its GAS for $F_0.$
\end{theorem}

For the sake of developing our theory in this paper, we find it beneficial to formalize a definition of enveloping in the following way: 

\begin{definition}\label{Def-EnvelopingV2} 
Let $\mathcal{I}$ be a nonempty interval and suppose the map $F_0$ in Eq. \eqref{Eq-F-Normalized} is from $ \mathcal{I}^k$ to $ \mathcal{I}.$ For any $x_1,\ldots,x_k\in \mathcal{I},$ define $\alpha:=\min\{x_1,\ldots,x_k\}$ and $ \beta:=\max\{x_1,\ldots,x_{k}\}.$ A continuous map $g\;:\;\mathcal{I}\to \mathcal{I}$ is called an enveloping of the map $F_0$  if (i) $g$ and $F_0$ share the fixed point $\bar y=1$ (ii) $g$ is decreasing  (iii) whenever $z=F_0(x_1,\ldots,x_k)$ exits the open interval $(\alpha,\beta),$ either $g(\beta)<z\leq \alpha$ or $g(\alpha)>z\geq \beta.$
\end{definition}
We considered the enveloping map $g$ decreasing, which is more stringent than the definition provided by El-Morshedy and Lopez \cite{El-Lo2008}, yet consistent with Cull's definition \cite{Cul2005} in one dimension.  Definition \ref{Def-EnvelopingV2} offers a benefit when enveloping a map and its expansions, as we clarify in the subsequent proposition: 
\begin{proposition}\label{Pro-hola-adios} 
Suppose $g$ envelopes $F$ in the equation $x_{n+1}=F(x_n,\ldots,x_{n-k+1}),$ where $F$ has a unique fixed point, then $g$ envelopes $\widetilde{F}$ in the equation 
$$y_{n+1}=\widetilde{F}(y_{n-j},\ldots,y_{n-j-k+1})=F_j(y_{n-j},\ldots,y_{n-j-k+1}).$$
\end{proposition}

\begin{proof}
We have $g,F$ and $\widetilde{F}$ share the same fixed point. Let $\beta=\max\{x_0,x_{-1},\ldots,x_{1-k}\}$ and $\alpha=\min\{x_0,x_{-1},\ldots,x_{1-k}\}.$ Since $g$ envelopes $F,$  we have  one of three compound inequalities 
$$g\left(\beta\right)<F(x_0,\ldots,x_{1-k})<\alpha\;\; \text{or}\;\; \alpha<F(x_0,\ldots,x_{1-k})<\beta\quad \text{or}\quad g\left(\alpha\right)>F(x_0,\ldots,x_{1-k})>\beta.$$ 
Now, $\widetilde{F}$ has $k$ components, but $k+j$ active variables. The middle inequality is the same for both $F$ and $\widetilde{F}.$ Now, we clarify the other two inequalities. Because $g$ is decreasing, we obtain 
$$g\left(\max\{x_0,\ldots,x_{1-k},\ldots,x_{1-k-j}\}\right)\leq g(\beta)\quad\text{and}\quad g\left(\min\{x_0,\ldots,x_{1-k},\ldots,x_{1-k-j}\}\right)\geq  g(\alpha).$$
Hence, $g$ envelopes $\widetilde{F}$.
\end{proof}

Next, we proceed to make the construction of $g$ more practical. 

\begin{definition}
 Let $g:\; \mathcal{I}\to \mathcal{I}$ be a map, and let $\Omega$ be a region in $\mathcal{I}\times \mathcal{I}.$  $g\in\Omega$ is used to denote that    $(x,g(x))\in\Omega$ for all $x\in\mathcal{I}.$
\end{definition}

\begin{lemma}\label{Lem-LocalDynamics}
Assume that $\bar x=1$ is a fixed point of $x_{n+1}=F_0(x_n,x_{n-1}).$  Let $(a_1,a_2)=\nabla F_0(1,1),$   and $M_1,M_2$ are the slopes of the curves $y=F_0(x,y)$ and $x=F_0(x,y)$ at $(1,1),$ respectively. If $\bar x=1$ is SLAS, then $|M_1|<1<|M_2|.$
\end{lemma}
\begin{proof}
The slopes of the curves of $y=F_0(x,y)$ and $x=F_0(x,y)$ at $(1,1)$ are given by
\begin{equation}\label{Eq-M1ANDM2}
    M_1=\frac{a_1}{1-a_2}\quad\text{and}\quad M_2=\frac{1-a_1}{a_2},
\end{equation}
respectively. Now, it is a matter of simple manipulations to find that  $|a_1|+|a_2|<1$ implies $|M_1|<1$ and $1<|M_2|.$
\end{proof}
Lemma \ref{Lem-LocalDynamics} tells us the local behavior of the curves of $y=F_0(x,y)$ and $x=F_0(x,y)$ near the LAS fixed point $\bar x=1.$ We aim to inflate the local dynamics under sufficient conditions to obtain global dynamics. Recall that if $\bar x=1$ is LAS under $F_0$, but not SLAS, then it must be SLAS under $F_j$ for some $j$ and, in this case, Lemma \ref{Lem-LocalDynamics} is valid on $F_j.$ Therefore, we will focus on the expansion $F_j$ that satisfies the LAS condition $|a_1|+|a_2|<1.$ Define the two-dimensional region 
\begin{equation}\label{Eq-RegionOfInterest}
\mathcal{R}:=\left\{(x,y):\; \alpha=\min\{x,y\}<F_j(x,y)<\max\{x,y\}=\beta\right\}\cup \{(1,1)\}.
\end{equation}
Since Lemma \ref{Lem-LocalDynamics} shows the slopes of the tangent lines at the point $(1,1),$ we conclude that $\mathcal{R}$ must have a positive area. Furthermore, we assume the curves of  $y=F(x,y)$ and $x=F(x,y)$ have a unique intersection in the interior of the positive quadrant.  This unique intersection must be at $x=y,$ and therefore, it must be  $(\bar x,\bar x)$.  Next, define the regions $\mathcal{R}_j,\; j=1,\ldots,4$ as follows: 
\begin{align*}
\mathcal{R}_1=&\{(x,y):\; F_j(x,y)<y\leq x\}\\
\mathcal{R}_2=&\{(x,y):\; F_j(x,y)<x\leq y\}\\
\mathcal{R}_3=&\{(x,y):\; F_j(x,y)>y\geq x\}\\
\mathcal{R}_4=&\{(x,y):\; F_j(x,y)>x\geq y\}.
\end{align*}
Fig. \ref{Fig-EnvelopingRegion} shows an illustration of the four regions together with the curves of $y=F_0(x,y)$ and $x=F_0(x,y).$
\definecolor{MyMaroon}{rgb}{128,0,0}
\definecolor{MyOrange}{rgb}{255,87,51}
\definecolor{ffvvqq}{rgb}{1,0.3333333333333333,0}
\definecolor{qqqqff}{rgb}{0.,0.,1.}
\definecolor{ffqqqq}{rgb}{1.,0.,0.}
\definecolor{cqcqcq}{rgb}{0.7529,0.7529,0.7529}
\definecolor{MyGreen}{rgb}{0,0.50196,0}
\definecolor{zzttqq}{rgb}{0.6,0.2,0}
\begin{figure}[htbp]
\centering
\begin{minipage}[t]{0.5\textwidth}
\begin{center}
\begin{tikzpicture}[line cap=round,line join=round,>=triangle 45,x=1.0cm,y=1.0cm,scale=1.8]
\draw[-triangle 45, line width=1.0pt,scale=1] (0,0) -- (3.3,0) node[below] {$x$};
\draw[line width=1.0pt,-triangle 45] (0,0) -- (-0.5,0);
\draw[-triangle 45, line width=1.0pt,scale=1] (0,0) -- (0,3.3) node[left] {$y$};
\draw[line width=1.0pt,-triangle 45] (0,0) -- (0.0,-0.5);
\fill[dash pattern=on 3pt off 3pt,color=ffqqqq,fill=ffqqqq,pattern=horizontal lines,pattern color=gray] (0.0,0.0) -- (0.023,0.174) -- (0.07,0.42) -- (0.135,0.633) -- (0.208,0.806) -- (0.2995,0.96) -- (0.4,1.084) -- (0.51,1.18) -- (0.61,1.256) -- (0.72,1.32) -- (0.854,1.37) -- (0.997,1.416) -- (1.103,1.442) -- (1.246,1.4686) -- (1.383,1.488) -- (1.5,1.5) -- (1.447,1.716) -- (1.376,2.014) -- (1.263,2.51) -- (1.204,2.78) -- (1.156,3.008) -- (1.1207,3.182) -- (1.079,3.3) -- (0,3.3) -- cycle;
\fill[dash pattern=on 3pt off 3pt,color=ffqqqq,fill=ffqqqq,pattern=vertical lines,pattern color=gray] (3.5,0) -- (1.912,0) -- (1.84,0.25) -- (1.77,0.4986) -- (1.696,0.754) -- (1.636,0.973) -- (1.59,1.152) -- (1.547,1.313) -- (1.5,1.5) -- (1.618,1.509) -- (1.7899,1.518) -- (1.999,1.524) -- (2.122,1.524) -- (2.276,1.523) -- (2.50,1.519) -- (2.626,1.515) -- (2.7976,1.509) -- (3.,1.5) -- (3.1535,1.49) -- (3.277,1.4857) -- (3.4984,1.473) -- cycle;
\draw[line width=1.0pt,color=blue,smooth] (0.0,0.0) -- (0.023,0.174) -- (0.07,0.42) -- (0.135,0.633) -- (0.208,0.806) -- (0.2995,0.96) -- (0.4,1.084) -- (0.51,1.18) -- (0.61,1.256) -- (0.72,1.32) -- (0.854,1.37) -- (0.997,1.416) -- (1.103,1.442) -- (1.246,1.4686) -- (1.383,1.488) -- (1.5,1.5) -- (1.618,1.509) -- (1.7899,1.518) -- (1.999,1.524) -- (2.122,1.524) -- (2.276,1.523) -- (2.50,1.519) -- (2.626,1.515) -- (2.7976,1.509) -- (3.,1.5) -- (3.1535,1.49) -- (3.277,1.4857) -- (3.4984,1.473);
\draw[line width=1.0pt,color=green,smooth] (0,0)--(3,3);
\draw[line width=1pt,color=magenta,dashed] (1.5,0)--(1.5,3.5);
\draw[line width=1pt,color=magenta,dashed] (0,3)--(3,0);
\draw[line width=1.0pt,color=red,smooth,samples=100,domain=1.08:1.95] plot(\x,{0-((2*(\x)^(2)+15*(\x)-36)/(2*(\x)+3))});
\draw[scale=1,rotate=0] (2.3,2.0) node[right] {\footnotesize $\mathcal{R}_1 $};
\draw[scale=1,rotate=0] (2,2.7) node[below] {\footnotesize $\mathcal{R}_2 $};
\draw[scale=1,rotate=0] (0.9,1.1) node[left] {\footnotesize $\mathcal{R}_3 $};
\draw[scale=1,rotate=0] (1,0.5) node[below] {\footnotesize $\mathcal{R}_4 $};
\draw[scale=1] (2.8,1.45) node[above] {\tiny  $y=F_j(x,y)$};
\draw[scale=1] (1.15,3.0) node[above,rotate=-77] {\tiny $x=F_j(x,y)$};
\draw[scale=1] (1.5,0.0) node[below] {\footnotesize $\bar x $};
\draw[scale=1] (3.0,0.0) node[below] {\footnotesize $2\bar x $};
\draw[scale=1] (0.0,3.0) node[left] {\footnotesize $2\bar x $};
\draw[line width=1.0pt,red,fill=yellow] (1.5,1.5) circle (1.5pt);
\draw[scale=1] (1.5,-0.5) node[below,rotate=0] {\footnotesize (i) $|M_1|<1$ and $M_2<-1$};
\end{tikzpicture}
\end{center}
\end{minipage}%
\begin{minipage}[t]{0.5\textwidth}
\begin{center}
\begin{tikzpicture}[line cap=round,line join=round,>=triangle 45,x=1.0cm,y=1.0cm,scale=1.8]
\draw[-triangle 45, line width=1.0pt,scale=1] (0,0) -- (3.5,0) node[below] {$x$};
\draw[line width=1.0pt,-triangle 45] (0,0) -- (-0.5,0);
\draw[-triangle 45, line width=1.0pt,scale=1] (0,0) -- (0,3.5) node[left] {$y$};
\draw[line width=1.0pt,-triangle 45] (0,0) -- (0.0,-0.5);
\draw[line width=1pt,color=magenta,dashed] (1.5,0)--(1.5,3.5);
\draw[line width=1pt,color=magenta,dashed] (0,3)--(3,0);
\fill[dash pattern=on 3pt off 3pt,color=ffqqqq,fill=ffqqqq,pattern=horizontal lines,pattern color=gray] (0,1.95) -- (0.114,1.9474) -- (0.245,1.938) -- (0.375,1.92) -- (0.489,1.902) -- (0.6045,1.8769) -- (0.718,1.847) -- (0.825,1.814) -- (0.924,1.779) -- (1.0165,1.743) -- (1.125,1.697) -- (1.228,1.648) -- (1.3296,1.596) -- (1.427,1.543) -- (1.5,1.5) -- (1.567,1.658) -- (1.627,1.81) -- (1.692,1.994) -- (1.745,2.1591) -- (1.7985,2.338) -- (1.8459,2.51) -- (1.8993,2.72) -- (1.941,2.897) -- (1.976,3.1) -- (0,3.1) -- cycle;
\fill[dash pattern=on 3pt off 3pt,color=ffqqqq,fill=ffqqqq,pattern=vertical lines,pattern color=gray] (0,0) -- (0.1452,0.066) -- (0.323,0.157) -- (0.47,0.244) -- (0.613,0.341) -- (0.755,0.453) -- (0.9096,0.596) -- (1.0637,0.7695) -- (1.188,0.936) -- (1.307,1.123) -- (1.41,1.318) -- (1.5,1.5) -- (1.614,1.43) -- (1.732,1.35) -- (1.85,1.268) -- (1.985,1.16) -- (2.12,1.053) -- (2.264,0.925) -- (2.371,0.83) -- (2.479,0.721) -- (2.59,0.607) -- (2.683,0.51) -- (2.776,0.4082) -- (2.86,0.314) -- (2.956,0.202) -- (3.035,0.108) -- (3.122498999199199,0) -- cycle;
\draw[line width=1.0pt,color=blue,smooth,samples=100,domain=0:3.2] plot(\x,{1.5*(13-3*\x*\x/(1.5^2))/10});
\draw[line width=1.0pt,color=red,smooth,samples=100,domain=0.0:2.0] plot(\x,{\x*(3*\x*\x/(1.5)^2+10)/(23-10*\x/1.5)});
\draw[scale=1,rotate=0] (2.3,2.0) node[right] {\footnotesize $\mathcal{R}_1 $};
\draw[scale=1,rotate=0] (2.2,2.6) node[below] {\footnotesize $\mathcal{R}_2 $};
\draw[scale=1,rotate=0] (0.9,1.1) node[left] {\footnotesize $\mathcal{R}_3 $};
\draw[scale=1,rotate=0] (0.8,0.75) node[below] {\footnotesize $\mathcal{R}_4 $};
\draw[scale=1,rotate=0] (2.7,0.45) node[above,rotate=-45] {\tiny  $y=F_j(x,y)$};
\draw[scale=1] (2.2,3.0) node[above,rotate=77] {\tiny $x=F_j(x,y)$};
\draw[line width=1.0pt,color=green,smooth] (0,0)--(3,3);
\draw[scale=1] (1.5,0.0) node[below] {\footnotesize $\bar x $};
\draw[scale=1] (3.0,0.0) node[below] {\footnotesize $2\bar x $};
\draw[scale=1] (0.0,3.0) node[left] {\footnotesize $2\bar x $};
\draw[line width=1.0pt,red,fill=yellow] (1.5,1.5) circle (1.5pt);
\draw[scale=1] (1.5,-0.5) node[below,rotate=0] {\footnotesize (ii) $|M_1|<1$ and $M_2>1$};
\end{tikzpicture}
\end{center}
\end{minipage}%
\caption{This figure illustrates the two scenarios of the curves $y=F_j(x,y)$ and $x=F_j(x,y)$. Both scenarios are possible, as we illustrate in our examples. $M_1$ and $M_2$ are defined in Eqs. \eqref{Eq-M1ANDM2}. The shaded region represents the region $\mathcal{R}$ as defined in Eq. \eqref{Eq-RegionOfInterest}. The dashed line $y=2\bar x-x$ is given to emphasize the importance of the slope $-1.$}\label{Fig-EnvelopingRegion}
\end{figure}
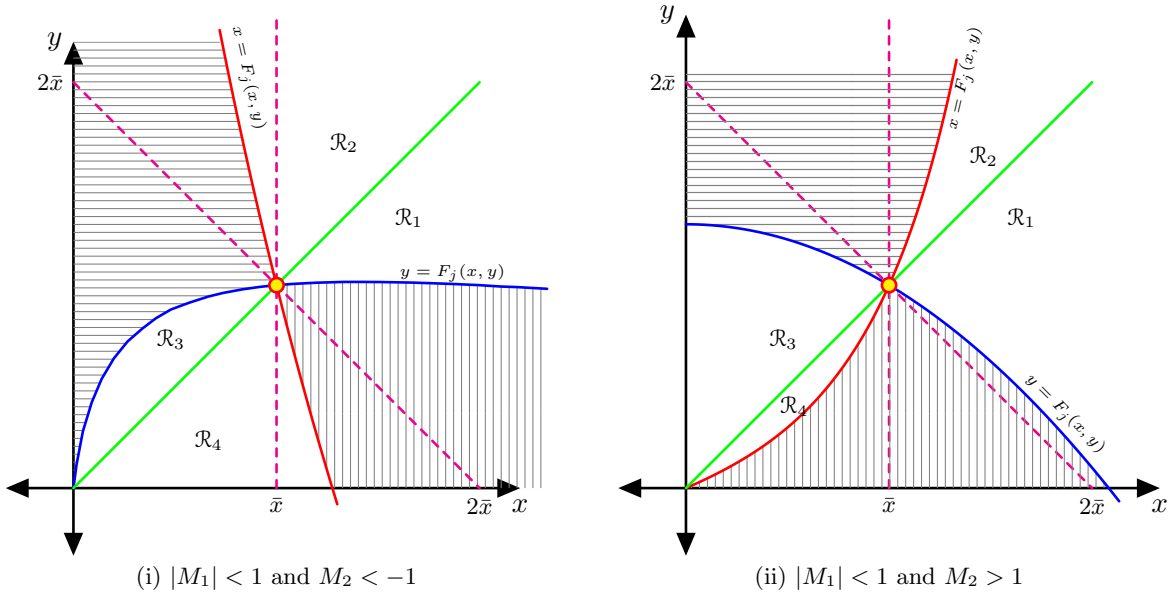

It is worth mentioning that it is possible for a one-dimensional map $g$ to envelop $F_0,$ but $g$ is not LAS at the fixed point. We clarify this in the following illustrative example: 

\begin{example}\label{Ex-DecreasingInTheSecond}\rm
Consider the difference equation 
$$x_{n+1}=F_0(x_n,x_{n-1})=\frac{a^2x_n}{(1+bx_n)(1+b(x_n+x_{n-1}))+abx_n},$$
where $b=\frac{1}{2}(a-1)$ and $a>1.$ It is straightforward to check that $\bar x=1$ is the only positive fixed point of $F_0.$ Furthermore,  because 
\begin{align*}
    |a_1|+|a_2|=&\|\nabla F_0(1,1)\|_1\\
    =& \frac{1}{4a^2}\left(|a^2-4a-1|+a^2-1\right)<1,
\end{align*}
then the fixed point is SLAS for all $a>1.$ In this example, we have 
$$M_1=\frac{a^2-4a-1}{1-5a^2}\quad\text{and}\quad M_2=-\frac{5a+1}{a+1}.$$
Thus, Part (i) of Fig. \ref{Fig-EnvelopingRegion} can be used for illustration.  
Now, define $g(x)=2-x$ when $0\leq x<2$ and zero when $x\geq 2.$ We have $g(1)=1,$ and every other point is either a $2$-cycle or eventually a $2$-cycle. So, the fixed point is not LAS for $g$. However, $g$ envelopes $F.$ Since conditions (i) and (ii) of Definition \ref{Def-EnvelopingV2} are clear, we proceed to clarify the third condition.
\begin{description}
\item{(i)} In $\mathcal{R}_1,$ we have $x\geq y,$ $\alpha=y$ and $\beta =x.$ Thus, we need to show that $g(x)=2-x<F_0(x,y).$  Observe that $F_0$ is decreasing in $y,$ and consequently, $F_0(x,y)\geq F_0(x,x)$ for all $x\geq y.$ Therefore, $F_0(x,y)+x\geq F_0(x,x)+x=:f_1(x).$ However, $f_1 $ is increasing in $x$ for all $x\geq 1.$ Because $f_1(x)>f_1(1)=2$ for all $x>1,$ then we obtain $2-x<F_0(x,y)$ in $\mathcal{R}_1.$ 
\item{(ii)} In $\mathcal{R}_2,$ we have $y\geq x,$ $\alpha=x,$  $\beta =y$ and $F_0(x,y)<x.$ Thus, we need to show that $g(y)=2-y<F_0(x,y).$ Alternatively, 
$$F_0(x,y)+y-2>0 \quad \Longleftrightarrow \quad p(x,y)=y^2+s_1(x)y-s_2(x)>0,$$
where
$$s_1(x)=\frac{(a-1)^2x^2+6(a-1)x-4(a-2)}{(a-1)((a-1)x+2)}\quad \text{and}\quad s_2(x)=\frac{2((a-1)^2x^2+2(a-2)x+4)}{(a-1)((a-1)x+2)}.$$
Also, the boundary curve $x=F_0(x,y)$ can be written as
$$h(x)=-x-\frac{2(a+1)}{a-1}+\frac{4a(a+1)}{(a-1)((a-1)x+2)}.$$
The proof's completion is a computational matter, and we omit the details. However, observe that $p(1,1)=0$ and $(1,1)$ is the only intersection between $p(x,y)=0$ and the boundary of $\mathcal{R}_2.$ Furthermore, $p(x,x)>0$ for all $x>1$ and $p(x,h(x))>0$ for all $x<1.$
\item{(iii)} In $\mathcal{R}_3,$ we have $y\geq x,$ $\alpha=x,$ $\beta =y$ and $F_0(x,y)>y.$ Thus, it remains to show that $g(x)=2-x>F_0(x,y).$ As in Part (i), since $F_0$ is decreasing in $y,$ we obtain $F_0(x,y)<F_0(x,x)$ for all $y>x.$ Thus, we proceed to show that $f_1(x)=F_0(x,x)+x<2$ for all $x<1.$ Indeed, $f_1$ is increasing in $x,$ and consequently, $f_1(x)<f_1(1)=2.$ This completes the task. 
\item{(iv)}  In $\mathcal{R}_4,$ we have $x\geq y,$ $\alpha=y,$ $\beta =x$ and $F_0(x,y)>x.$ Thus, it remains to show that $g(y)=2-y>F_0(x,y).$ This is similar to case (ii), but with switched inequalities. We skip the computations because they are similar.  
\end{description}
\end{example}
Observe that the domain of the map $F_0$ in Example \ref{Ex-DecreasingInTheSecond} is $\mathbb{R}_+^2$; however, it is straightforward to obtain 
$$F_0(x,y)\leq F_0(x,0)\leq F_0\left(\frac{2}{a-1},0\right).$$
So, we can ignore the first three iterates and consider $F_0:\mathcal{I}^2\to\mathcal{I},$ where $\mathcal{I}:=\left[0,F_0\left(\frac{2}{a-1},0\right)\right].$

When the fixed point is SLAS, the curve of the enveloping map that serves our interest must be in the closure of the shaded regions of Fig. \ref{Fig-EnvelopingRegion} as we prove in the following result:

\begin{proposition}\label{Pr-NecessaryCondForEnv}
Consider $F_j$ to be an expansion of $F_0$ that satisfies the LAS condition $\|\nabla F_j(1,1)\|_1<1.$  If $y=g(x)$ envelopes $F_j,$ and $g$ has no $2$-cycle, then we must have
$$\min\{x,g(x)\}< F_j(x,g(x))< \max\{x,g(x)\}\quad\text{for all}\quad x\neq 1.$$
\end{proposition}
\begin{proof}
The LAS condition $\|\nabla F_j(1,1)\|_1<1$ implies that we have one of the scenarios in Fig. \ref{Fig-EnvelopingRegion}.  The conclusion means the curve of the function $y=g(x)$ is contained in the region $\mathcal{R}.$ So, assume to the contrary that an enveloping function $y=g(x)$ is not completely inside the region $\mathcal{R}.$ If part of $y=g(x)$ intersects the region $\mathcal{R}_1 ,$ then there exists a point $(x,y)\in\mathcal{R}_1$ such that $y<g(x).$ This contradicts the fact that we must have $g(x)<F_j(x,y)<y<x$.  Similarly, if part of the curve $y=g(x)$ intersects the region $\mathcal{R}_3$.  Since $g$ is decreasing, $y=g(x)$ cannot intersect regions $\mathcal{R}_2$ and $\mathcal{R}_4$ of Case (ii) in Fig. \ref{Fig-EnvelopingRegion}. Thus, it remains to focus on regions $\mathcal{R}_2$ and $\mathcal{R}_4$ of Case (i) in Fig. \ref{Fig-EnvelopingRegion}.  Suppose that part of the curve of  $y=g(x)$ intersects the region $\mathcal{R}_2,$ then there exists a point $(x,g(x))\in\mathcal{R}_2$ such that  $F_j(x,g(x))<x<g(x).$ Since $g$ envelopes $F_j,$ we obtain 
$$g^2(x)<F_j(x,g(x))<x<g(x).$$
Therefore, we obtain two subsequences that converge to a $2$-cycle of $g,$ i.e., $g^{2n}(x)\searrow \bar y$ and $g^{2n+1}(x)\nearrow g(\bar y),$ where
$\{\bar y,g(\bar y)\}$ is a $2$-cycle of $g$ and 
$$\bar y\leq x<\bar x<g(\bar y).$$
This violates our assumption. The argument is similar if $y=g(x)$ intersects  the region $\mathcal{R}_4,$ which completes the proof.  
\end{proof}
Proposition \ref{Pr-NecessaryCondForEnv} provides a useful geometric insight for constructing an enveloping map $g$, prompting us to explore this idea further in the sequel.  When a two-dimensional map $F$ is increasing in its second argument, i.e., $F(\cdot,\uparrow),$ and satisfies the LAS condition, Lemma \ref{Lem-LocalDynamics} implies $M_2>1.$ However, the monotonic nature in the second argument here makes LAS and SLAS equivalent, as we show in the following lemma:  
\begin{lemma}\label{Lem-LAS=SLAS}
Let $\bar x$ be a fixed point of a differentiable two-dimensional map $F.$ Suppose $F$ is increasing in its second argument in a neighborhood containing $\bar x.$ Then $\bar x$ is LAS if and only if it is SLAS. 
\end{lemma}
\begin{proof} 
SLAS implies that LAS is a consequence of Theorem \ref{Th-Al-Ca2024}. However,
let $(a_1,a_2)=\nabla F(\bar x,\bar x), $ where $a_2>0.$  The eigenvalues of the Jacobian matrix are the zeros of $p(\lambda)=\lambda^2-a_1\lambda-a_2,$ which are real and in the interval $(-1,1)$ when $p(-1)>0$ and $p(1)>0.$ Now, since $a_2>0,$ we have
$$p(-1)>0,\;p(1)>0\quad \Longleftrightarrow\quad  |a_1|+|a_2|=\|\nabla F\|_1<1,$$
which completes the proof. 

 \end{proof}
 
 Since $M_2=\frac{1-a_1}{a_2}>1$ for maps of the form $F(\cdot,\uparrow),$ this scenario  fits under Part (ii) of Fig. \ref{Fig-EnvelopingRegion}, and it can be used for illustration of the following result: 
\begin{theorem} \label{Th-MT1}
Let $F_j:\;\mathcal{I}\times\mathcal{I}\to \mathcal{I}$ be an expansion of $F_0$ such that $F_j$ is increasing in its second argument. Suppose the fixed point $\bar x=1$ is unique and LAS. If $g:\;\mathcal{I}\to\mathcal{I}$ is a continuous decreasing function such that $\min\{x,g(x)\}< F_j(x,g(x))< \max\{x,g(x)\}$ for all $x\neq 1,$ then $g$ is an enveloping of $F_j.$
\end{theorem}
\begin{proof}
As clarified in Proposition \ref{Pr-NecessaryCondForEnv}, we obtain $(x,g(x))\in \mathcal{R}$ for all $x.$ Also, each of the systems $(y,x)=(g(x),F_j(x,y))$ and $(y,y)=(g(x),F_j(x,y))$ has a unique solution, namely $x=y=1$. Now, we proceed to verify condition (iii) of Definition \ref{Def-EnvelopingV2}. We split the proof into cases based on the regions $\mathcal{R}_1$ to $\mathcal{R}_4,$ but first, we define the regions 
$$\mathcal{R}_l=\{(x,y)\in\mathcal{R}:\; x<1\}\quad\text{and}\quad 
\mathcal{R}_r=\{(x,y)\in\mathcal{R}:\; x>1\}.$$
\begin{description}
\item{(i)} Let $(x,y)\in\mathcal{R}_1.$ We have $x>y>F_j(x,y)$ and we need to show that $F_j(x,y)>g(x).$ In this case, $(x,g(x))\in \mathcal{R}_r,$ and consequently, $g(x)<F_j(x,g(x)) .$  Also, in $\mathcal{R}_1,$ we have $g(x)<y.$ Because $F_j$ is increasing in its second argument, we obtain 
$F_j(x,y)>F_j(x,g(x)),$ which completes the proof of this case. 
\item{(ii)} Let $(x,y)\in\mathcal{R}_2.$ We have $y>x>F_j(x,y)$ and we need to show that $F_j(x,y)>g(y).$ In this case, $(y,g(y))\in \mathcal{R}_r,$ and consequently, $g(y)<F_j(x,g(y)) .$ Since $g(y)<y$ and $F_j$ is increasing in its second argument, we obtain $F_j(x,g(y))<F_j(x,y),$ which completes the proof of this case. 
\item{(iii)} Let $(x,y)\in\mathcal{R}_3.$ We have $x<y<F_j(x,y)$ and we need to show that $g(x)>F_j(x,y).$ In this case, $(x,g(x))\in \mathcal{R}_l,$ and consequently, $F_j(x,g(x))<g(x) .$ Also, we have  $g(x)>y,$ and the monotonicity of $F_j$ in its second argument gives us   $F_j(x,g(x))>F_j(x,y),$ which completes the proof of this case. 
\item{(iv)} Let $(x,y)\in\mathcal{R}_4.$ We have $y<x<F_j(x,y)$ and we need to show that $F_j(x,y)<g(y).$ In this case, $g(y)>y$ and the monotonicity of $F_j$ in its second argument gives us   $F_j(x,g(y))>F_j(x,y).$ Finally, the point $(x,g(y))$ must be above the curve $y=F_j(x,y),$ which implies $F_j(x,g(y))<g(y).$ This completes the proof. 
\end{description}
\end{proof}
Theorem \ref{Th-MT1} tells us that we can obtain an enveloping of $F_j$ by constructing a continuous decreasing map $g$ such that $g\in\mathcal{R}.$ Based on this fact and Theorem \ref{Th-Lopez}, we obtain the following corollary.
\begin{corollary}
Consider the assumptions of Theorem \ref{Th-MT1}. If there exists a continuous decreasing function $g:\;\mathcal{I}\to\mathcal{I}$ such that $g\in \mathcal{R},$ $\bar x=1$ is LAS for $g$ and $g$ has no $2$-cycle, then $\bar x=1$ is globally attracting for $F_0.$
\end{corollary}
\begin{proof}
From Proposition \ref{Pr-2CycleGAS}, $\bar x$ is GAS for $g.$ Since $g$ is an enveloping of $F_j$ by Theorem \ref{Th-MT1}, $\bar x$ is GAS for $F_j$ by Theorem \ref{Th-Lopez}. Therefore, $\bar x$ is GAS for $F_0.$
\end{proof}

\section{Enveloping  in contrast to embedding}\label{Sec-Enveloping}
When a map $F$ is monotonic in each argument, we say it is of mixed monotonicity. The dynamics of this category of maps have been investigated extensively in the literature \cite{Al2022,Al-Ca-Ka2024B}, and the embedding technique has been developed to establish global stability. The idea hinges on finding an appropriate partial order that can be used to extend the map $F$ into a higher-dimensional map that is increasing with respect to the defined partial order. Even though this notion has been extensively investigated, exploring it within the context of the enveloping technique is intriguing. In this section, we use $\uparrow$ and $\downarrow$ to denote strictly increasing and strictly decreasing, respectively.  We begin by defining the partial order $\leq_\tau$ that is compatible with the monotonicity of $F$ \cite{Al-Ca-Ka2024B}, and then, we extract the main theorem. We define
$$(x_1,x_2,\ldots,x_k)\leq_\tau (u_1,u_2,\ldots,u_k)$$
if and only if $x_i\leq u_i$ ($u_i\leq x_i$) whenever $F$ is increasing (decreasing) in its $i^{th}$ argument. Since our focus in the sequel will be limited to $2$-dimensions, we need to formalize some concepts in the following definition \cite{Al-Ca-Ka2024B}:
\begin{definition}\label{Def-Pseudo}
For a two-dimensional function $F$ that is monotonic in each of its arguments, define the point $ P_\tau = (u_1, u_2),$ where $ u_i = x $ if $F$ is increasing in its \( i^{th} \) argument, and \( u_i = y \) otherwise. $P_\tau^t$ denotes the ``dual" of $P_\tau$, whereby $x$ and $y$ are interchanged. When $x \leq y$, the points $P_\tau$ and $P_\tau^t$ denote the extreme points of a boxed region in a $2$-dimensional space under the $\tau$-partial order.  A solution $(x,y),\; x\neq y$ of the system $(x,y)=(F(P_\tau), F(P^t_\tau))$ is called a pseudo-fixed point of $F.$ 
\end{definition}

\begin{theorem}\cite{Al-Ca-Ka2024B} \label{Th-GS-embedding}
Let $F:\;\mathcal{I}^2\to \mathcal{I}$ be continuous and monotonic in each of its arguments. Suppose $F$ has a unique fixed point and no pseudo-fixed points. If for each initial condition $X_0$ in the domain of $F,$ there exists a point $P_\tau(x,y)$ such that 
\begin{equation}\label{Eq-embedding1}
x\leq y,\;\; x<F(P_\tau), \;\; y>F(P^t_\tau) \quad \text{and}\quad P_\tau\leq_\tau X_0 \leq_\tau P_\tau^t,
\end{equation}
then the fixed point of $F$ is a global attractor. 
\end{theorem}
It is interesting to observe that Theorem \ref{Th-GS-embedding} may fail because of the vacuous solution of the inequalities in \eqref{Eq-embedding1}. A simple example that clarifies this note is the map $F(x,y)=\frac{bx}{1+y},\; b>1$. This serves as our rationale for examining maps of mixed monotonicity in two dimensions and contrasting the embedding and enveloping methods.

 Now, we are in a position to focus on each of the four cases: $F(\uparrow,\uparrow),$ $F(\downarrow,\uparrow),$ $F(\downarrow,\downarrow)$ and $F(\uparrow,\downarrow).$ Keep in mind that LAS and SLAS are equivalent in the first two cases (Lemma \ref{Lem-LAS=SLAS}).
\subsection{The case {$\mathbf{F(\uparrow,\uparrow)}$}}\label{Sec-SubSection4.1}
The standard practice in the literature is to use the one-dimensional map $g(x)=F(x,x)$ as the mean to investigate the dynamics of $F.$ We stress that we assume a unique positive fixed point at $\bar x=1,$ and $\bar x$ is LAS under $F.$ Unlike the case in Eq. \eqref{Eq-Linear1}, the SLAS of $\bar x$ is a fact that follows from Lemma \ref{Lem-LAS=SLAS}.
In other words, looking for an expansion of $F_j$ in this case is unnecessary. Also, 
this case is within the scope of both Theorem \ref{Th-GS-embedding} and Theorem \ref{Th-MT1}.
To invoke Theorem \ref{Th-GS-embedding}, we explicitly give the partial order
$$(x,y)\leq_\tau (u,v)\quad \iff x\leq u\quad \text{and}\quad y\leq v.$$
The points are $P_\tau=(x,x)$ and $P_\tau^t=(y,y),$ where $x\leq y.$ Pseuduo-fixed points are not possible in this case because the only solution of $(x,y)=(F(P_\tau),F(P_\tau^t))$ is $x=y=1.$ Therefore, all that is needed is to investigate the inequalities in \eqref{Eq-embedding1} of Theorem \ref{Th-GS-embedding}. The following proposition simplifies the technicalities of applying theorems \ref{Th-GS-embedding} and \ref{Th-MT1}.

\begin{proposition}\label{Pr-IncrIncr}
Let $\Omega_1=\{(x,y):\quad x<1,y>1\},$ $\Omega_2=\{(x,y):\quad x>1,y<1\},$ $\Omega_3=\{(x,y):\quad x<y,\; x<F(P_\tau),F(P_\tau^t)<y\},$ and $\mathcal{R}$ as defined in Eq. \eqref{Eq-RegionOfInterest}. Then $\Omega_1\subset \Omega_3,$ $\Omega_1\cup\Omega_2\subset \mathcal{R}$ and $\Omega_3\subset \mathcal{R}.$
\end{proposition}
\begin{proof}
Because for $x<1,$ we have $F(x,x)<F(1,1)=1,$ and for $y>1,$ we have $F(y,y)>1,$ then the three inequalities that define $\Omega_3$ are valid. This shows that $\Omega_1\subset \Omega_3.$ Next, we show that $\Omega_1\cup\Omega_2\subset \mathcal{R}.$
From the uniqueness of the fixed point and its local stability, we must have $F(x,x)>x$ for $x<1$ and $F(x,x)<x$ for $x>1.$ Now, in $\Omega_1,$ we have
$$x<F(x,x)\leq  F(x,y)\leq F(y,y)<y.$$
Therefore, $\Omega_1\subset \mathcal{R}.$ Similarly, in $\Omega_2,$ we have
$$y<F(y,y)\leq  F(x,y)\leq F(x,x)<x,$$
and consequently, $\Omega_2\subset \mathcal{R}.$

Finally, we show that $\Omega_3\subset \mathcal{R}.$ Let $(x,y)\in\Omega_3.$ This and the nature of monotonicity in $F$  give us 
$$x<F(x,x)<F(x,y)<F(y,y)<y,$$
which completes the proof. 
\end{proof}

From the fact that $\Omega_1\subset \Omega_3,$ for each initial condition $X_0,$ we can pick $P_\tau$ so that $P_\tau\leq_\tau X_0\leq P^t_\tau.  $ Hence, LAS implies GAS based on Theorem \ref{Th-GS-embedding}.
Next, we turn our attention to the enveloping technique. 
Proposition \ref{Pr-IncrIncr} and Theorem \ref{Th-MT1} tell us that we can define a generic continuous decreasing map $g:\; [0,\infty)\to [0,\infty)$ in the region $\Omega_1\cup\Omega_2\cup\{(1,1)\},$ which serves as an enveloping map for $F.$ For instance, both of the following maps achieve the task: 
\begin{equation}\label{Eq-Phi1Phi2}
    \phi_1(x)=\frac{a(1+x)}{1+(2a-1)x},\quad \phi_2(x)=\frac{a+1}{ax+1},\quad a>1.
\end{equation}
Note that in comparison to $\phi_1$ and $\phi_2$ in Eq. \eqref{Eq-OneDimensional-g1}, we fixed $b_1=a,c_1=2a-1,$ $c_2=a$ and $b_2=0.$
Furthermore, $\bar x=1$ is globally attracting under both of them. As an example of this case, consider 
$$x_{n+1}=F(x_n,x_{n-1})=\frac{a(bx_n+1)(bx_{n-1}+1)}{b(bx_n+2)(bx_{n-1}+2)},$$
where $a>0$ and $b$ is the unique positive solution of the equation 
$$x^3+(4-a)x^2+2(2-a)x-a=0.$$
It is a matter of simple computations to check that $\bar x=1$ is the unique positive fixed point, and the enveloping technique also shows that $\bar x=1$ is GAS under $F.$
\subsection{The case $\mathbf{F(\downarrow,\uparrow)}$}
Again, this case is within the scope of Theorem \ref{Th-GS-embedding} and
Theorem \ref{Th-MT1}. We proceed under the assumption that $\bar x=1$ is the unique positive fixed point. Also, LAS and SLAS are equivalent based on Lemma \ref{Lem-LAS=SLAS}. The partial order in Theorem \ref{Th-GS-embedding} is defined by
$$(x,y)\leq_\tau (u,v)\quad \iff u\leq x\quad \text{and}\quad y\leq v.$$
The point $P_\tau=(y,x),$ where $x<y.$ Thus, it remains to explore the solution region of the inequalities 
$x<F(y,x)$ and $F(x,y)<y.$ 

\begin{proposition}
Let $F(\downarrow,\uparrow),$ and suppose that $\bar x=1$ is unique and LAS. The region obtained by the inequalities $x<y,\; x<F(y,x)$ and $F(x,y)<y$ is contained in the region $\mathcal{R}.$ 
\end{proposition}
\begin{proof}
The inequality $F(x,y)<y$ is used to define $\mathcal{R}$ when $x<y.$ Also, from the monotonicities in $F,$ we have $F(y,x)<F(x,y).$ This completes the proof. 
\end{proof}
Thus, Theorem \ref{Th-GS-embedding} is applicable here, and there is no need to seek an expansion of $F.$ Also, Theorem \ref{Th-MT1} can be applied by seeking an enveloping map $g$ that fits within the shaded region of Part (i) of Fig. \ref{Fig-EnvelopingRegion}. In fact, the map $g(x)=F(x,x)$ can be a good candidate when the decreasing argument of $F$ dominates the increasing one. We avoid giving the details because they follow the ideas of Subsection \ref{Sec-SubSection4.1}. However, recall that Lemma \ref{Lem-LAS=SLAS} shows that $\bar x=1$ is LAS iff it is SLAS. Also,  since $F$ is increasing in the second argument, we have 
$$\|\nabla F(1,1)\|_1=|F_x(1,1)|+F_y(1,1)<1\quad \text{implies}\quad |g^\prime(1)|<1.$$
The following example clarifies the case of this subsection.  
\begin{example} \rm
Consider the second-order difference equation 
$$x_{n+1}=F_0(x_n,x_{n-1})=\frac{\frac{a}{2}(x_{n-1}+1)}{1+x_{n-1}+(a-2)x_{n}^2},\quad a>2.$$
We have $F_0:\;\mathbb{R}_+^2\to \mathbb{R}_+$ and $F_0(\downarrow,\uparrow).$ Furthermore, $F_0$  has a unique positive fixed point for all $a>2$, namely $\bar x=1.$ $\bar x=1$ is SLAS when 
$$\|\nabla F_0(1,1)\|_1=\frac{2(a-2)}{a}+\frac{a-2}{2a}<1.$$
This is satisfied as long as $2<a<\frac{10}{3}.$ Lemma \ref{Lem-LAS=SLAS} shows that SLAS is equivalent to LAS here. Indeed,
the Jacobian matrix at $\bar x=1$ is given by 
$$J=\frac{1}{2a}\left[
      \begin{array}{cc}
        4(2-a) &a-2\\
        2a&0\\
      \end{array}
    \right].$$
    It is straightforward to check that the eigenvalues of $J$ are inside the open unit disk as long as $2<a<\frac{10}{3}.$  Thus, $F_0$ satisfies the LAS condition, and this example fits case (ii) of Fig. \ref{Fig-EnvelopingRegion}.  Theorem \ref{Th-MT1} can be applied to obtain GAS through enveloping, and Theorem \ref{Th-GS-embedding} can be applied to obtain GAS through embedding.   
\end{example}

\subsection{The case $\mathbf{F(\downarrow,\downarrow)}$}
In this case, the one-dimensional map $g(x)=F(x,x)$  is useless without assuming that $F$ satisfies the LAS condition. Modifying the example in Eq. \eqref{Eq-Linear1} to read as 
\begin{equation}\label{Eq-Linear2}
    x_{n+1}=F(x_n,x_{n-1})=-\frac{3}{5}x_n-\frac{3}{5}x_{n-1}
\end{equation}
will clarify this argument. 
The zero fixed point is LAS for $F,$ but unstable for $g(x)=F(x,x).$ However, we explore this case within the context of embedding and enveloping techniques. The partial order in Theorem \ref{Th-GS-embedding} is given by
$$(x,y)\leq_\tau (u,v)\quad \iff\quad u\leq x\quad \text{and}\quad v\leq y.$$
The point $P_\tau=(y,y)$ and $P_\tau^t=(x,x),$ where $x\leq y.$ Pseuduo-fixed points are  solutions of $(x,y)=(F(P_\tau),F(P_\tau^t))$ for which $x\neq y.$ The following proposition addresses this issue:
\begin{proposition}\label{Pr-DecrDecr1}
Assume $F$ satisfies the LAS condition at the unique positive fixed point, and define $g(x)=F(x,x).$ Then $g$ has no $2$-cycle iff $F$ has no pseudo-fixed points.  
\end{proposition}
\begin{proof}
Since $g$ has no $2$-cycle, then the only solution of $(x,y)=(g(y),g(x))$ is the fixed point of $F.$ Therefore, the system $(x,y)=(F(P_\tau),F(P_\tau^t))$ has no  solution other than $x=y=1.$
\end{proof}
Next, we clarify the relationship between the solution of the inequalities in \eqref{Eq-embedding1} of Theorem \ref{Th-GS-embedding} and the region $\mathcal{R}$ that we heavily use in the enveloping technique.  

\begin{proposition}\label{Pr-DecrDecr2}
Let $\Omega_i,i=1,2,3$ be as defined in Proposition \ref{Pr-IncrIncr}
 and $\mathcal{R}$ as defined in Eq. \eqref{Eq-RegionOfInterest}. Then 
 $\Omega_3\subset \mathcal{R},$ $\mathcal{R} \subset\Omega_1\cup\Omega_2\cup\{(1,1)\}$ and $\Omega_3\subset \Omega_1.$  
\end{proposition}
\begin{proof}
Let $(x,y)\in \Omega_3.$ We have $x<y,x<F(y,y)$ and $y>F(x,x).$ These inequalities and the fact that $F$ is decreasing in both arguments give us
$$x<F(y,y)<F(x,y)<F(x,x)<y.$$
Therefore, $\Omega_3\subset \mathcal{R}.$ Next, let $(x,y)\in \mathcal{R}.$ The case $x=y$ gives us $(x,y)=(1,1).$ So, we proceed with $x<y.$ This means $x<F(x,y)<y.$ Because $F$ is decreasing in both arguments and by implicit differentiation of the curves $x=F(x,y)$ and $y=F(x,y),$ we find the slope of the two curves is negative. This means the two curves are bound to be in $\Omega_1.$ Similarly, if $x>y,$ then $y<F(x,y)<x,$ and we obtain $(x,y)\in\Omega_2.$  Finally, $\Omega_3\subset \Omega_1$ becomes clear. 
\end{proof}

If $F$ fails to satisfy the LAS condition at $\bar x$, then we can take its first expansion, i.e., $F_1(x,y)=F(F(x,y),x).$ In this case, $F_1$ is increasing in its second argument, and we can treat $F_1$ based on the machinery of the previous two subsections.   As an illustrative example of this case, consider 
\begin{equation}\label{Eq-DecDec}
x_{n+1}=F(x_n,x_{n-1})=\frac{(b+1)^2}{(bx_n+1)(bx_{n-1}+1)},
\end{equation}
where $b>0.$  $\bar x=1$ is a fixed point, and since $f(x)=\frac{(b+1)^2}{(bx+1)^2}$ is decreasing,  $\bar x=1$ must be the unique positive fixed point. Because the Jacobian matrix at $\bar x=1$ is
$$
J:=\begin{bmatrix}
-\frac{b}{b+1}&-\frac{b}{b+1}\\
1&0\\
\end{bmatrix},
$$
$F$ satisfies the LAS condition for all $0<b<1.$ However, $\bar x=1$ is LAS for all $b>0.$ Interestingly, the fixed point is LAS under the one-dimensional map $g(x)=F(x,x)$ as long as $0<b<1.$ At $b=1,$ a $2$-cycle is born, and consequently, $g$ fails to help us detect stability. On the other hand, Proposition \ref{Pr-DecrDecr1} tells us that pseudo-fixed points of $F$ are created at $b=1,$ which makes Theorem \ref{Th-GS-embedding} inapplicable.

As mentioned earlier, to stay armed with the LAS condition, we invoke the first expansion of $F,$ which is
\begin{equation}
\begin{split}
y_{n+1}=&F_1(y_{n-1},y_{n-2})\\
=&\frac{(b+1)^2(by_{n-2}+1)}{b(b+1)^2+(by_{n-1}+1)(by_{n-2}+1)},\quad y_{-2},y_{-1},y_0\;\text{and}\;b>0.
\end{split}
\end{equation}
However, the nature of monotonicity in $F_1$ differs from $F.$ Indeed, we have $F_1(\downarrow,\uparrow).$  Thus, it becomes under the jurisdiction of Theorem \ref{Th-MT1} after taking Proposition \ref{Pro-hola-adios} into account. It is a matter of simple computations to show that the graph of $g(x)=\frac{b+1}{bx+1}$ belongs to the region $\mathcal{R},$ and consequently, $g$  achieves the enveloping task, which assures the GAS of $\bar x=1$ in Eq. \eqref{Eq-DecDec}.

\subsection{The case $\mathbf{F(\uparrow,\downarrow)}$}
In this case, the partial order that makes Theorem \ref{Th-GS-embedding} applicable is as follows:
$$(x,y)\leq_\tau (u,v)\quad \iff\quad x\leq u\quad \text{and}\quad v\leq y.$$
The point $P_\tau=(x,y),$ where $x<y.$ The next result establishes a connection between the solution of the inequalities in Theorem \ref{Th-GS-embedding} and the region $\mathcal{R}$ in Eq. \eqref{Eq-RegionOfInterest}.

\begin{proposition}\label{Pr-IncrDecr}
Let $\Omega:=\{(x,y):\; x<y,\; x<F(P_\tau),\; y>F(P_\tau^t)\},$ and consider $\mathcal{R}$ as defined in Eq. \eqref{Eq-RegionOfInterest}. Then $\Omega\subset \mathcal{R}.$ 
\end{proposition}
\begin{proof}
If $\Omega$ is empty, then there is nothing to show. So, let $(x,y)\in\Omega.$ Since $x<y,$ we need to show that $x<F(x,y)<y.$ But this is straightforward because $x<F(x,y)$ is already satisfied, and the type of monotonicity in the map $F$ gives us 
$$F(x,y)<F(y,y)<F(y,x)<y.$$
This completes the proof.
\end{proof}
\begin{lemma}\label{Lem-IncrDecr}
Consider $F$ to be ${F(\uparrow,\downarrow)}$ and $\mathcal{R}$ as defined in Eq. \eqref{Eq-RegionOfInterest}. Let $\phi:\;\mathcal{I}\to\mathcal{I}$ be decreasing and differentiable map such that $(x,\phi(x))\in\mathcal{R}$ for all $x\in \mathcal{I}$. If $(x,\phi^{-1}(x))\in\mathcal{R},$ then $g(x)=F(\phi(x),x)$ is an enveloping of $F.$
\end{lemma}
\begin{proof}
Since $\phi\in \mathcal{R}$ and it is decreasing, we have $\phi(1)=1,$ $\phi^{-1}$ is well defined, $g$ is decreasing and $g(1)=1.$ Therefore, it remains to check condition (iii) of Definition \ref{Def-EnvelopingV2}. As in the proof of Theorem \ref{Th-MT1}, we split the proof based on the regions $\mathcal{R}_i.$ 
\begin{description}
\item{(i)} Let $(x,y)\in\mathcal{R}_1.$ We have $x>y>F(x,y)$, and we need to show that $F(x,y)>g(x).$ Because $x>y>\phi(x),$ we obtain 
$$F(x,y)>F(x,x)>F(\phi(x),x)=g(x).$$ 
\item{(ii)} Let $(x,y)\in\mathcal{R}_2.$ We have $y>x>F(x,y)$ and we need to show that $F(x,y)>g(y).$ In this case, we have $y>\phi^{-1}(x)$ and $x>\phi(y).$ The nature of monotonicity in $F$ gives us
$$F(x,y)>F(\phi(y),y)=g(y).$$
\item{(iii)} Let $(x,y)\in\mathcal{R}_3.$ We have $x<y<F(x,y)$ and we need to show that $F(x,y)<g(x).$ Since $x<1,$ we have $x<\phi(x),$ and based on the monotonicities in $F,$ we obtain 
$$F(x,y)<F(x,x)<F(\phi(x),x)=g(x).$$ 
\item{(iv)} Let $(x,y)\in\mathcal{R}_4.$ We have $y<x<F(x,y)$ and we need to show that $F(x,y)<g(y).$ Because $y<x<1,$ we have $x<\phi(x)<\phi(y).$ Therefore, we obtain 
$$F(x,y)<F(\phi(x),y)<F(\phi(y),y))=g(y).$$
\end{description}
\end{proof}

\begin{proposition}
Suppose $F$ satisfies the LAS condition at the unique fixed point $\bar x=1,$ and assume the map $y=\phi(x)$ in Lemma \ref{Lem-IncrDecr} is obtained from the equation $y=F(y,x).$ Then $\bar x$ is LAS for $\phi$ and the following statements are equivalent.
\begin{description}
\item{(a)} $F$ has a pseudo-fixed point
\item{(b)} $\phi$ has a $2$-cycle
\item{(c)} $g(x)=F(\phi(x),x)$ has a $2$-cycle.
\end{description}
\end{proposition}
\begin{proof}
As in Lemma \ref{Lem-LocalDynamics}, implicit differentiation of $y=F(y,x)$ gives us 
$$|\phi^\prime(1)|=\left|\frac{a_2}{1-a_1}\right|,$$
which is smaller than one due to the SLAS. Now, we show the equivalency between the three statements. 
\begin{description}
\item{(a) $\Longrightarrow$ (b)} A pseudo-fixed point of $F$ means $(x,y)=(F(x,y),F(y,x))$ has a solution for some $x\neq y.$ Observe that each equation is an inverse relation of the other, and since $\phi$ is obtained from $\phi(x)=F(\phi(x),x),$ $\{x,y\}$ is a $2$-cycle of $\phi.$ 
\item{(b) $\Longrightarrow$ (c)} Let $C_2:=\{x_0,x_1\}$ be a $2$-cycle of $\phi.$ This means $\phi(x_0)=x_1$ and $\phi(x_1)=x_0,$ where $x_0\neq x_1.$ From the structure of $\phi,$ we have \begin{align*}
\phi(x_0)=&F(\phi(x_0),x_0)=g(x_0)=F(x_1,x_0)=x_1\\
\phi(x_1)=&F(\phi(x_1),x_1)=g(x_1)=F(x_0,x_1)=x_0.
\end{align*}
Therefore, $C_2$ is a $2$-cycle of $g.$
\item{(c) $\Longrightarrow$ (a)} The proof is similar and omitted. 
\end{description}
\end{proof}
\begin{theorem}\label{Th-IncrDecr}
Assume $\bar x=1$ is a unique fixed point of $F$, and it is SLAS. Consider $\mathcal{R}$ as defined in Eq. \eqref{Eq-RegionOfInterest}, and let $\phi:\;\mathcal{I}\to\mathcal{I}$ be decreasing and differentiable map such that $(x,\phi(x))\in\mathcal{R}$ for all $x$. Suppose $\bar x=1$ is LAS under $\phi$ and $(x,\phi^{-1}(x))\in\mathcal{R}.$ If $g(x)=F(\phi(x),x)$ has no $2$-cycle, then $\bar x$ is GAS for $F.$
\end{theorem}
\begin{proof}
Since $\bar x=1$ is SLAS under $F$ and LAS under $\phi,$ we obtain 
$$|g^\prime(1)|=|F_x(1,1)\phi^\prime(1)+F_y(1,1)|\leq \|\nabla F(1,1)\|_1<1.$$
So, the fixed point is LAS under $g.$ Because $g$ has no $2$-cycle, then $\bar x$ is GAS under $g,$ and because $g$ envelopes $F,$ then the fixed point is GAS under $F.$ 
\end{proof}
Next, we provide our final example, which clarifies the case in this subsection.   

\begin{example}\rm
The Ricker map with delay and stocking 
$$x_{n+1}=F(x_n,x_{n-1})=x_n\exp(b-x_{n-1})+h,\;b,h>0$$
was considered in \cite{Al-Ka2025}. $F$ has a unique positive fixed point that we denote by $\bar x$ ($\bar x>h$).  To be consistent, we normalize the fixed point and write
$$y_{n+1}=F_0(y_n,y_{n-1})=y_n\exp{(b-\bar x y_{n-1})}+\frac{h}{\bar x},\;b,h>0.$$
We obtain SLAS at  $\bar y=1$ when  
$$\|\nabla F_0(1,1)\|_1=\left(1-\frac{h}{\bar x}\right)+(\bar x-h)=\left(\bar x-h\right)\left(1+\frac{1}{\bar x}\right)<1.$$
This implies $\bar y$ is SLAS as long as $\bar x<\frac{1}{2}(h+\sqrt{h^2+4h}):=h^*.$ This is the same GAS condition obtained in \cite{Al-Ka2025} based on Theorem \ref{Th-GS-embedding}. Next, we discuss the enveloping technique by applying Lemma \ref{Lem-IncrDecr} and Theorem \ref{Th-IncrDecr}.  Consider the map 
$$y=\phi(x)=\frac{h}{\bar x(1-\exp{(b-\bar x x)})}.$$
Note that we can consider $F_0:\; (0,\infty)\to(0,\infty).$ Also, we obtained the map $\phi$ from $y=F_0(y,x),$ which is used in the boundary of the region $\Omega$ in Proposition \ref{Pr-IncrDecr}. Furthermore, the SLAS of $\bar y$ forces $|\phi^\prime(1)|<1$. On the other hand, $x=F_0(x,y)$, which also serves to define $\Omega$, represents $\phi^{-1}$. Because Proposition \ref{Pr-IncrDecr} shows that $\phi$ and $\phi^{-1}$ are both in $\mathcal{R}.$ Lemma \ref{Lem-IncrDecr} tells us that  
$$g(x)=F_0(\phi(x),x)=\frac{h\exp{(b-\bar x x)}}{\bar x(1-\exp{(b-\bar x x)})}+h$$
envelopes $F_0.$ It is a computational matter to  show that $g$ has no $2$-cycle as long as $\bar x<h^*.$  
\end{example}

\section{Conclusion}

The enveloping technique is built on finding a one-dimensional map that can bind the orbits of a multidimensional map so that convergence in the one-dimensional map forces convergence in the multidimensional map. The nonexistence of a $2$-cycle in the one-dimensional map plays a crucial role in the global stability. However, the existing literature results lack the practical apparatus to achieve the goal.

In this paper, we develop theoretical and geometric machinery for the two-dimensional case, building upon the expansion strategy presented in \cite{Al-Ca2024}. When a hyperbolic fixed point of a map $F_0$ is locally asymptotically stable, there is an expansion $F_j$ of $F_0$ such that $\|\nabla F_j\|_1<1$ at the fixed point. In this case, we call the fixed point (Definition \ref{Def-LASCondition}) strongly locally asymptotically stable (SLAS) under $F_j$, and the key idea is to use $F_j$ rather than $F_0$ in tackling global stability. The relationship between $\|\nabla F_j\|_1<1$ and LAS in $F_0$  was elucidated in reference \cite{Al-Ca2024}.  Also, the illustrated cases presented in Section \ref{Sec-Enveloping} clarify the connection. However, as a simple example, consider $F_0(x,y)=-\frac{3}{5}x-\frac{3}{5}y.$ In this instance, the zero fixed point is LAS and $\|\nabla F_0\|_1>1$. But, $F_1(x,y)=-\frac{6}{25}x+\frac{9}{25}y$, resulting in $\|\nabla F_1\|_1=\frac{3}{5}<1.$

Our enveloping approach leads to a practical geometric approach (as summarized in Fig. \ref{Fig-EnvelopingRegion}, Theorem \ref{Th-MT1} and Theorem \ref{Th-IncrDecr}) that helps us in addressing global stability in the two-dimensional case. We connected the enveloping and embedding techniques within the context of SLAS. We thoroughly characterized global stability in maps of mixed monotonicity and maps that are increasing in the second argument.

This paper contains numerous examples that illustrate all aspects of our developed results. We were not interested in the complete analysis of the dynamics of the considered examples but in illustrating the various aspects of our developed theory. 
\medskip

\noindent{\textbf{Acknowledgement:}} We are grateful to the anonymous reviewers for their insightful comments and recommendations, which improved a previous version of this manuscript. Also, the authors thank Victor Jimenez Lopez for our stimulating discussion about the dominance technique. The first author is supported by sabbatical leave from the American University of Sharjah and a Maria Zambrano grant for attracting international talent from the Polytechnic University of Cartagena. 
\bibliographystyle{unsrt}

\bibliography{AlSharawi-bibliography2}
\end{document}